\documentclass{article}
\usepackage{amsmath}
\usepackage{amsfonts}
\usepackage{graphics}
\usepackage{mathrsfs}
\usepackage{latexsym}
\usepackage{flafter}
\usepackage{amssymb}
\usepackage{dsfont}
\usepackage{colortbl}
\usepackage{multirow}
\usepackage{tabularx}
\usepackage{cite}
\allowdisplaybreaks[3]
\newtheorem{theorem}{Theorem}

\begin{document}
\title{\textbf{A periodic FM-BEM for solving the acoustic transmission problems in periodic media}}
\author{Wenhui Meng\thanks{E-mail address: mwh@nwu.edu.cn(Wenhui Meng).}, ~Ruifei Liu\\
\small\parbox{120mm}{\begin{center}\emph{School of Mathematics, Northwest University, Xi'an, 710127, China}
\end{center}}}
\date{}
\maketitle

\begin{abstract}
This paper presents a new fast multipole boundary element method (FM-BEM) for solving the acoustic transmission problems in 2D periodic media. We divide the periodic media into many fundamental blocks, and then construct the boundary integral equations in the fundamental block. The fast multipole algorithm is proposed for the square and hexagon periodic systems, the convergence of the algorithm is analyzed. We then apply the proposed method to the acoustic transmission problems for liquid phononic crystals and derive the acoustic band gaps of the phononic crystals. By comparing the results with those from plane wave expansion method, we conclude that our method is efficient and accurate.

\textbf{keywords}: fast multipole method, boundary element method, acoustic transmission problem, phononic crystal, acoustic band gap
\end{abstract}

\section{Introduction}
Acoustic transmission problem is widely used in various areas, such as sonar, phononic crystal and nondestructive testing. Acoustic problem is a classical problem that was frequently solved by the boundary element method (BEM). Kress \cite{1,2,3,4} proposed an indirect BEM for solving acoustic scattering and transmission problems, the solution of the problem is expressed as the form of single and/or double layer potentials, boundary integral equation is then constructed by the boundary condition. However, since BEM produces dense and nonsymmetric linear system of equation, its computational efficiency has been a serious problem for solving large-scale models.

The fast multipole method (FMM) \cite{5,6} is an effective method to accelerate the solution of BEM. In recent years, the fast multipole boundary element method (FM-BEM) has been developed to solve a variety of problems that are computationally intensive, such as potential problems, Stokes flow problems and acoustic wave problems. Some applications of FM-BEM in acoustic wave problems can be found in \cite{7,8,9,10}. In \cite{7}, Liu summarized the theory of FM-BEM and its applications in engineering, this book is particularly useful to researchers newly working on this subject.

Wave propagation in periodic media is frequently encountered in engineering problems, such as phononic crystals and photonic crystals. Several methods have been used to study the wave transmission characteristics in periodic media, such as the plane wave expansion (PWE) method \cite{11}, transfer matrix (TM) method \cite{12}, finite difference time domain (FDTD) method \cite{13}, multiple-scattering theory (MST) \cite{14}, and so on. But these methods have certain drawbacks and limitations. Compared to the above methods, FM-BEM has the characteristics of fast solving speed, high precision and wide adaptability.

In recent years, BEM has been widely used to analyze the acoustic and electromagnetic problems in periodic medium \cite{15,16,17,18,19,20}. However, there is little existing work studying the applications of FMM in this field. Otani and Nishimura devote themselves to using the FM-BEM to solve acoustic and electromagnetic scattering problems in periodic medium \cite{18,19,20}. Since the transmission problem leads to mixed boundary conditions, it follows that the boundary integral equation is more complex than that of the scattering problem.

In this paper, a periodic FM-BEM is proposed for solving the acoustic transmission problem in the infinite periodic medium. In Section 2, according to the periodic arrangement of scatterers, the periodic media is divided into many equal fundamental blocks, each fundamental block covers $M$ scatterers. The transmission problem for infinite number of scatterers is then transformed into that for $M$ scatterers in the fundamental block. For the transmission boundary condition, we construct the boundary integral equations by the method of Kress. Section 3 present a periodic FMM for solving the boundary integral equations, and analyze the error of the algorithm. Finally, the present method is applied to calculate the acoustic band gaps of water and mercury phononic crystals, exact band gaps are obtained.

\section{The periodic boundary integral equations}
Consider the time-harmonic acoustic wave propagation in a composite media with periodic fiber arrangements of an infinite spatial length. The matrix and fibers are homogeneous and isotropic, the cross section of fibers may wish to be denoted by the region $\Omega_i(i=1,2,\cdots)$ and the boundary curve of $\Omega_i$ is $\Gamma_i$. Suppose the length of the periodic media is $L$ in the x-direction and is infinite in the y-direction (see Fig.1).
\begin{figure}[htbp]
\centering
\scalebox{0.9}{\includegraphics{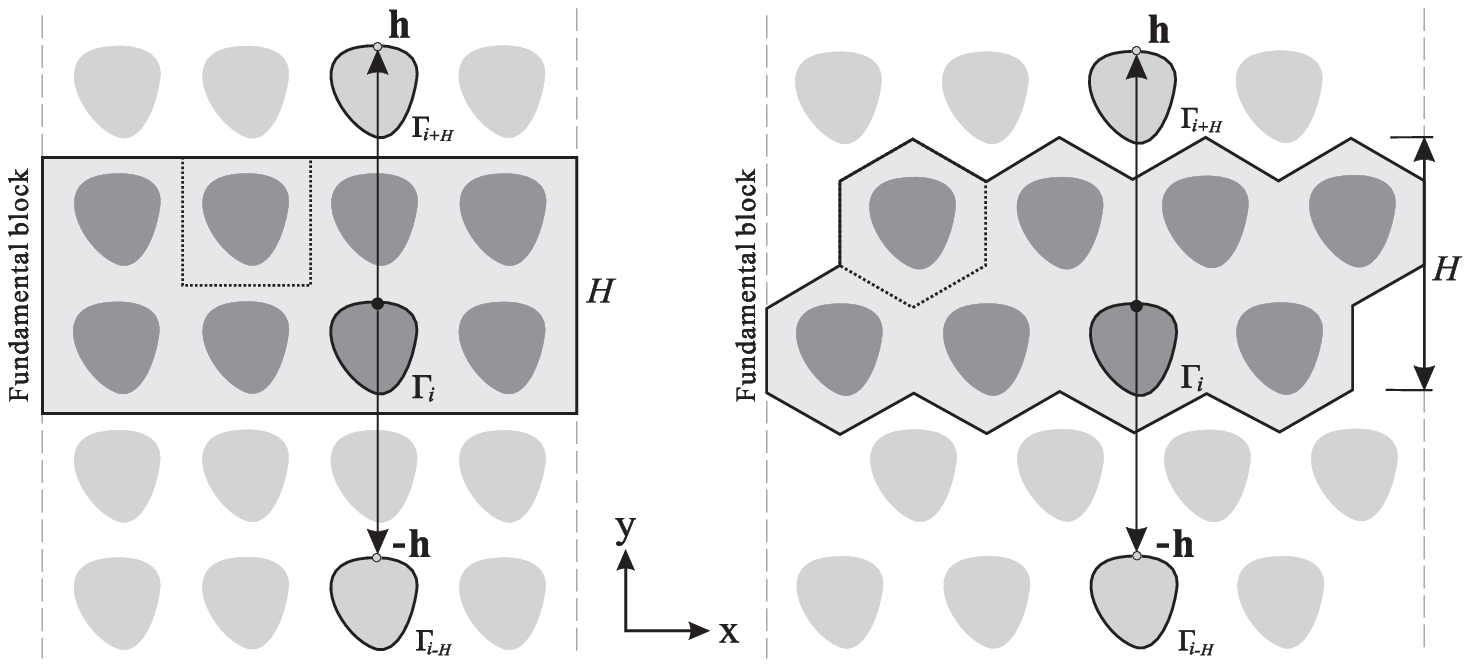}}\\
\small Fig.1 The square periodic structure (left) and the hexagon periodic structure (right).
\end{figure}

Consider the incident wave $u^{inc}(\mathbf{x})=e^{ik\mathbf{x}\cdot{\mathbf{d}}}$ with $\mathbf{d}=(\cos\theta,\sin\theta)$, where $\theta$ is the incident angle. Thus, the refracted wave $u_i(\mathbf{x})$ and scattered wave $u^s(\mathbf{x})$ obey the following Helmholtz equations:
\begin{eqnarray}
\left\{\begin{array}{ll}
\Delta{u^s}(\mathbf{x})+k^2u^s(\mathbf{x})=0,& ~~\mathbf{x}\in{\Omega_0},\\
\Delta{u_i}(\mathbf{x})+k_1^2u_i(\mathbf{x})=0,& ~~\mathbf{x}\in{\Omega_i},i=1,2,\cdots,
\end{array}\right.
\end{eqnarray}
where $\Omega_0={\mathds{R}^2}\backslash\bigcup\overline{\Omega_i}$ denotes the region of the matrix. Suppose $c_1$ and $c_2$ are wave velocities in scatterer and matrix respectively, thus we have $k_1=\omega/c_1$ and $k=\omega/c_2$, where $\omega$ is the radian frequency of the incident wave.

The continuity of the pressure and of the normal velocity across the interface leads to the transmission boundary conditions
\begin{equation}
u_i(\mathbf{x})=u(\mathbf{x}),\qquad\frac{1}{\rho_1}\frac{\partial{u}_i}{\partial\nu}(\mathbf{x})=\frac{1}{\rho_2}\frac{\partial{u}}{\partial\nu}(\mathbf{x}),\qquad \mathbf{x}\in\Gamma_i,i=1,2,\cdots,
\end{equation}
where $\nu$ is the unit outward normal to $\Gamma_i$, $u(\mathbf{x})$ is the total field
$$u(\mathbf{x})=u^{inc}(\mathbf{x})+u^s(\mathbf{x}),\qquad \mathbf{x}\in{\Omega_0},$$
and $\rho_1, \rho_2$ are the densities of scatterer and matrix respectively.

Divide the periodic media into many equal fundamental blocks with length $H$ in the y-direction. In each fundamental block, total of $M$ fibers are arranged to form a lattice, containing $M_H$ level rows and $M_L$ vertical columns of fibers, such as Fig.1 for $M_H=2,M_L=4$. This leads to the following quasi-period conditions:
\begin{equation}
u(\mathbf{x}\pm\mathbf{h})=\alpha{u}(\mathbf{x}),\qquad\frac{\partial{u}}{\partial\nu}(\mathbf{x}\pm\mathbf{h})=\alpha\frac{\partial{u}}{\partial\nu}(\mathbf{x}),
\qquad\mathbf{x}\in\Gamma_i,i=1,2,\cdots,M,
\end{equation}
where $\alpha=e^{ik\mathbf{h}\cdot{\mathbf{d}}}$ and $\mathbf{h}=(0,H)$.

To solve the problem (1)-(3), first introduce the single and double-layer potentials:
$$S_l\varphi_i(\mathbf{x}):=\int_{\Gamma_i}\Phi_l(\mathbf{x},\mathbf{y})\varphi_i(\mathbf{y})ds(\mathbf{y}),\qquad \mathbf{x}\in{\mathds{R}^2\backslash{\Gamma_i}},$$
$$D_l\phi_i(\mathbf{x}):=\int_{\Gamma_i}\frac{\partial\Phi_l(\mathbf{x},\mathbf{y})}{\partial\nu(\mathbf{y})}\phi_i(\mathbf{y})ds(\mathbf{y}),
\qquad\mathbf{x}\in{\mathds{R}^2\backslash{\Gamma_i}},$$
where $\varphi_i$ and $\phi_i$ are continuous functions defined on $\Gamma_i$, $\Phi_1$ and $\Phi_2$ are the fundamental solutions of Helmholtz equation, expressed as
$$\Phi_1(\mathbf{x},\mathbf{y})=\frac{\mathrm{i}}{4}H^{(1)}_0(k_1|\mathbf{x}-\mathbf{y}|),\qquad\mathbf{x}\neq{\mathbf{y}}.$$
$$\Phi_2(\mathbf{x},\mathbf{y})=\frac{\mathrm{i}}{4}H^{(1)}_0(k|\mathbf{x}-\mathbf{y}|),\qquad \mathbf{x}\neq{\mathbf{y}}.$$
We seek the solution to (1)-(3) in the form of combined single and double-layer potentials:
\begin{eqnarray}
&&u^s(\mathbf{x})=\displaystyle{\sum_{i=1}^\infty\Big\{D_2\phi_i(\mathbf{x})-\mathrm{i}\eta{S}_2\varphi_i(\mathbf{x})\Big\}},\qquad{\mathbf{x}}\in{\Omega_0},\\
&&u_i(\mathbf{x})=D_1\phi_i(\mathbf{x})-\mathrm{i}\eta{S}_1\varphi_i(\mathbf{x}),\qquad\mathbf{x}\in{\Omega_i},i=1,2,\cdots,
\end{eqnarray}
where $\eta$ is a real constant. Let the solutions $u^s(\mathbf{x})$ and $u_i(\mathbf{x})$ satisfy the boundary conditions (2), according to the jump relations for single and double-layer potentials\cite{1}, we derive the following integral equations:
\begin{eqnarray}
\left\{\begin{array}{l}
\displaystyle\phi_i(\mathbf{x})-{D}_1\phi_i(\mathbf{x})+\mathrm{i}\eta{S}_1\varphi_i(\mathbf{x})+\sum_{j=1}^\infty\Big\{{D}_2\phi_j(\mathbf{x})
-\mathrm{i}\eta{S}_2\varphi_j(\mathbf{x})\Big\}=-u^{inc}(\mathbf{x}),\\
\displaystyle\frac{\mathrm{i}\eta(1+\varrho)}{2}\varphi_i(\mathbf{x})-\varrho{T}_1\phi_i(\mathbf{x})+\mathrm{i}\varrho\eta{K}_1\varphi_i(\mathbf{x})+\sum_{j=1}^\infty
\Big\{{T}_2\phi_j(\mathbf{x})-\mathrm{i}\eta{K}_2\varphi_j(\mathbf{x})\Big\}=-\frac{\partial{u}^{inc}}{\partial\nu}(\mathbf{x}),
\end{array}\right.
\end{eqnarray}
where $\mathbf{x}\in\Gamma_i(i=1,2,\cdots)$, $\varrho=\rho_2/\rho_1$, and
$${K}_l\varphi_i(\mathbf{x}):=\int_{\Gamma_i}\frac{\partial\Phi_l(\mathbf{x},\mathbf{y})}{\partial\nu(\mathbf{x})}\varphi_i(\mathbf{y})ds(\mathbf{y}),\qquad
\mathbf{x}\in{\mathds{R}^2\backslash{\Gamma_i}},$$
$${T}_l\phi_i(\mathbf{x}):=\frac{\partial}{\partial\nu(\mathbf{x})}\int_{\Gamma_i}\frac{\partial\Phi_l(\mathbf{x},\mathbf{y})}{\partial\nu(\mathbf{y})}
\phi_i(\mathbf{y})ds(\mathbf{y}),\qquad \mathbf{x}\in{\mathds{R}^2\backslash{\Gamma_i}}.$$
in addition, by the period conditions (3), the scattered field $u^s(\mathbf{x})$ can be written as the form:
\begin{eqnarray}
u^s(\mathbf{x})\hspace{-0.6cm}&&=\sum_{i=1}^\infty\int_{\Gamma_i}\left\{\frac{\partial\Phi_2(\mathbf{x},\mathbf{y})}
{\partial\nu(\mathbf{y})}\phi_{i}(\mathbf{y})-\mathrm{i}\eta\Phi_2(\mathbf{x},\mathbf{y})\varphi_{i}(\mathbf{y})\right\}ds(\mathbf{y})\nonumber\\
&&=\sum_{i=1}^M\sum_{m=-\infty}^\infty\alpha^m\int_{\Gamma_i}\left\{\frac{\partial\Phi_2(\mathbf{x},\mathbf{y}+m\mathbf{h})}
{\partial\nu(\mathbf{y})}\phi_i(\mathbf{y})-\mathrm{i}\eta\Phi_2(\mathbf{x},\mathbf{y}+m\mathbf{h})\varphi_i(\mathbf{y})\right\}ds(\mathbf{y}).
\end{eqnarray}
Thus, for $\mathbf{x}\in\Gamma_i(i=1,2,\cdots,M)$, the integral equations (6) can be converted into
\begin{eqnarray}
\left\{\begin{array}{l}
\displaystyle\phi_i(\mathbf{x})-{D}_1\phi_i(\mathbf{x})+\mathrm{i}\eta{S}_1\varphi_i(\mathbf{x})+\sum_{j=1}^M\Big\{\overline{D}_2\phi_j(\mathbf{x})-\mathrm{i}\eta
\overline{S}_2\varphi_j(\mathbf{x})\Big\}=-u^{inc}(\mathbf{x}),\\
\displaystyle\frac{\mathrm{i}\eta(1+\varrho)}{2}\varphi_i(\mathbf{x})-\varrho{T}_1\phi_i(\mathbf{x})+\mathrm{i}\varrho\eta{K}_1\varphi_i(\mathbf{x})+\sum_{j=1}^M
\Big\{\overline{T}_2\phi_j(\mathbf{x})-\mathrm{i}\eta\overline{K}_2\varphi_j(\mathbf{x})\Big\}=-\frac{\partial{u}^{inc}}{\partial\nu}(\mathbf{x}),
\end{array}\right.
\end{eqnarray}
where $\overline{S}_2,\overline{D}_2,\overline{K}_2,\overline{T}_2$ are the periodic potentials:
\begin{eqnarray*}
&&{\overline{S}}_2\varphi_i(\mathbf{x}):=\sum_{m=-\infty}^{\infty}\alpha^m\int_{\Gamma_i}\Phi_2(\mathbf{x},\mathbf{y}+m\mathbf{h})\varphi_i(\mathbf{y})ds(\mathbf{y}),\\
&&{\overline{D}}_2\phi_i(\mathbf{x}):=\sum_{m=-\infty}^{\infty}\alpha^m\int_{\Gamma_i}\frac{\partial\Phi_2(\mathbf{x},\mathbf{y}+m\mathbf{h})}
{\partial\nu(\mathbf{y})}\phi_i(\mathbf{y})ds(\mathbf{y}),\\
&&{\overline{K}}_2\varphi_i(\mathbf{x}):=\sum_{m=-\infty}^{\infty}\alpha^m\int_{\Gamma_i}\frac{\partial\Phi_2(\mathbf{x},\mathbf{y}+m\mathbf{h})}
{\partial\nu(\mathbf{x})}\varphi_i(\mathbf{y})ds(\mathbf{y}),\\
&&{\overline{T}}_2\phi_i(\mathbf{\mathbf{x}}):=\sum_{m=-\infty}^{\infty}\alpha^m\frac{\partial}{\partial\nu(\mathbf{x})}\int_{\Gamma_i}\frac{\partial\Phi_2(\mathbf{x},
\mathbf{y}+m\mathbf{h})}{\partial\nu(\mathbf{y})}\phi_i(\mathbf{y})ds(\mathbf{y}).
\end{eqnarray*}
In \cite{22}, we prove that $\overline{S}_2,\overline{D}_2,\overline{K}_2,\overline{T}_2$ are uniformly convergent when $kH(1\pm\sin\theta)\neq2n\pi(n\in\mathds{N})$, but the convergence rates are $\mathcal{O}(1/\sqrt{q})$, where $q$ is the truncation order. We may adjust to improve the convergence rate of the scheme.

For $\mathbf{x}\in\Gamma_i(i=1,2,\cdots,M)$, the boundary integral equations (8) can be written as
\begin{eqnarray}
(I-A)\Psi=B,
\end{eqnarray}
where
$$A=\left(\begin{array}{ccccccc}
D_1-\overline{D}_2 & \mathrm{i}\eta(\overline{S}_2-S_1) & -\overline{D}_2 & \mathrm{i}\eta\overline{S}_2 & \cdots & \cdots \\
\displaystyle\frac{2\mathrm{i}(\overline{T}_2-\varrho{T}_1)}{\eta(1+\varrho)} & \displaystyle\frac{2(\overline{K}_2-\varrho{K}_1)}{1+\varrho} & \displaystyle\frac{2\mathrm{i}\overline{T}_2}{\eta(1+\varrho)} & \displaystyle\frac{2\overline{K}_2}{1+\varrho} & \cdots & \cdots \\
\vdots & \vdots & \ddots & \ddots & \vdots & \vdots \\
\cdots & \cdots & -\overline{D}_2 & \mathrm{i}\eta\overline{S}_2 & D_1-\overline{D}_2 & \mathrm{i}\eta(\overline{S}_2-S_1) \\
\cdots & \cdots & \displaystyle\frac{2\mathrm{i}\overline{T}_2}{\eta(1+\varrho)} & \displaystyle\frac{2\overline{K}_2}{1+\varrho} & \displaystyle\frac{2\mathrm{i}(\overline{T}_2-\varrho{T}_1)}{\eta(1+\varrho)} & \displaystyle\frac{2(\overline{K}_2-\varrho{K}_1)}{1+\varrho}
\end{array}
\right)$$
and
$$\Psi=(\phi_1,\varphi_1,\phi_2,\varphi_2,\cdots,\phi_M,\varphi_M)^\top,$$
$$B=(-u^{inc}\Big|_{\Gamma_1},\frac{2\mathrm{i}\partial{u}^{inc}}{\eta(1+\varrho)\partial\nu}\Big|_{\Gamma_1},\cdots,-u^{inc}\Big|_{\Gamma_M},
\frac{2\mathrm{i}\partial{u}^{inc}}{\eta(1+\varrho)\partial\nu}\Big|_{\Gamma_M})^\top.$$

It can be seen that (9) is a Fredholm system of integral equations of the second kind, thus it is uniquely solvable. In next section, a periodic FMM is proposed to accelerate the solution of the equations.

\section{The periodic FMM}

\subsection{The periodic tree structure}
Section 3.1 shows that we can only construct the tree structure in the fundamental block and its adjacent blocks. Since the scatterers are periodic arrangement in the fundamental block, it follows that we can construct a tree structure in a unit lattice and then copy it to other unit lattice of the fundamental block. We call the tree structure in unit lattice as the basic tree (see Fig.2).

For the square lattice, the traditional square quad-tree structure is available. In the basic tree, the relationships between cells can be determined by traditional way. But how to determine the relationships between the cells in different basic tree?
\begin{figure}[htbp]
\centering
\scalebox{0.65}{\includegraphics{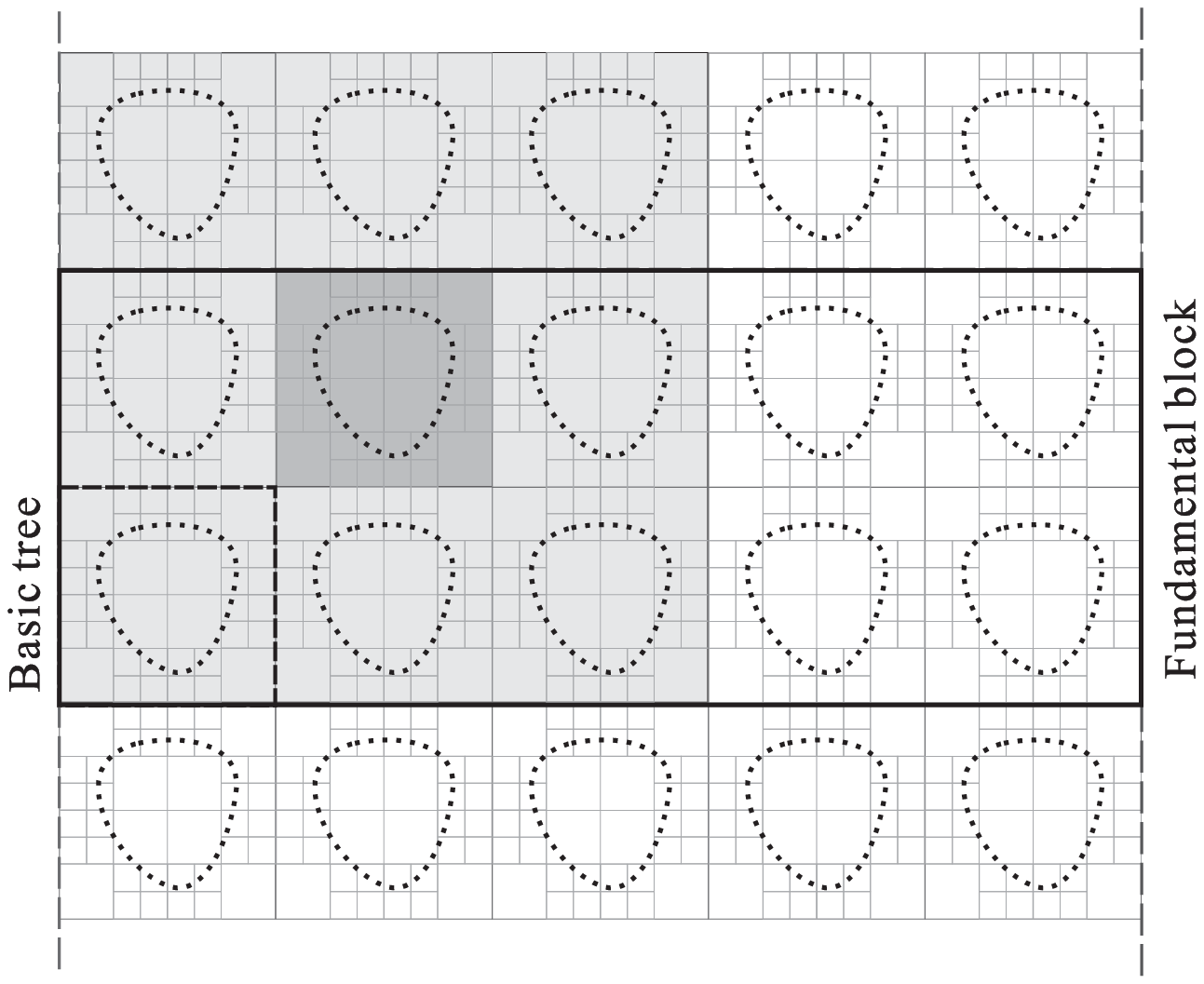}}\\
\small Fig.2 The periodic tree structure for the square unit lattice.
\end{figure}

Suppose that $Q_s$ and $Q_f$ are the basic trees including the source point $\mathbf{x}$ and field point $\mathbf{y}$ respectively. If $Q_s$ and $Q_f$ share at least one vertex, then they are called adjacent trees, otherwise called far trees.

(1) If $Q_s$ and $Q_f$ are adjacent, then the relationships of their cells can be determined by traditional way. Fig.4 shows the M2M, M2L and L2L translations between two adjacent basic trees.

\begin{figure}[htbp]
\centering
\scalebox{0.85}{\includegraphics{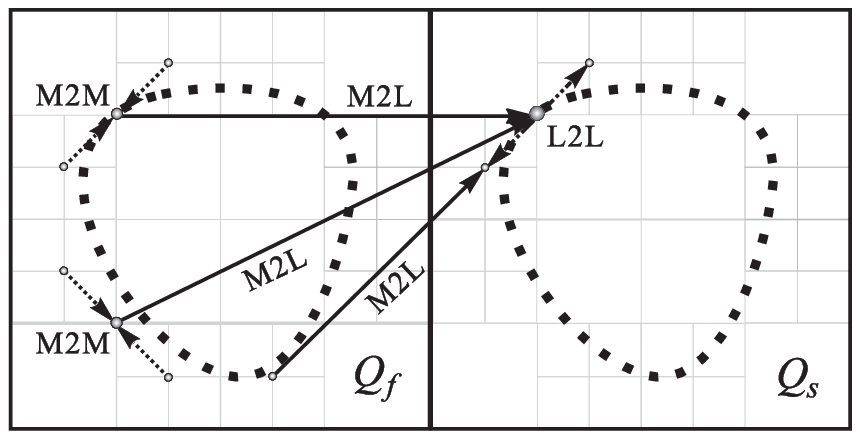}}\\
\small Fig.3 FMM translations between the adjacent basic trees.
\end{figure}

(2) If $Q_s$ and $Q_f$ are far, then their level 0 cells are well separated or far away. The M2L translation can be applied between the level 0 cells. As shown in Fig.4.
\begin{figure}[htbp]
\centering
\scalebox{0.8}{\includegraphics{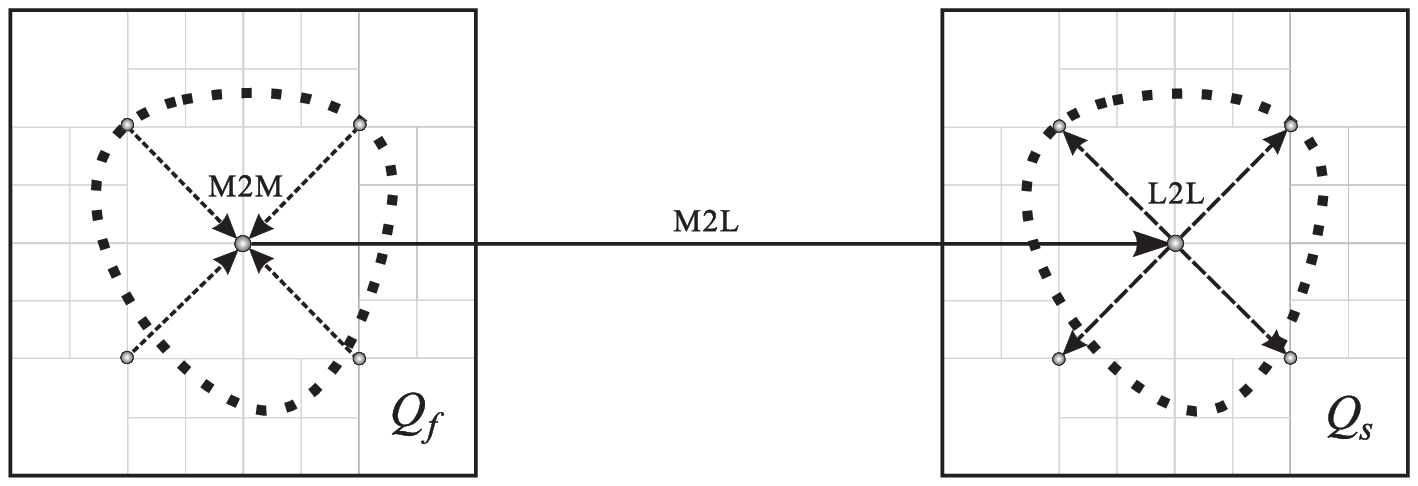}}\\
\small Fig.4 FMM translations between the far basic trees.
\end{figure}

For the hexagon lattice, we can construct regular triangle quad-tree structure. We first divide the regular hexagonal lattice (level 0 cell) into six equal regular triangles (level 1 cells), and then divide each regular triangle (level 1 cell) into four small equal regular triangles (level 2 cells). In this way, a quad-tree structure can also be constructed (see Fig.5).
\begin{figure}[htbp]
\centering
\scalebox{0.85}{\includegraphics{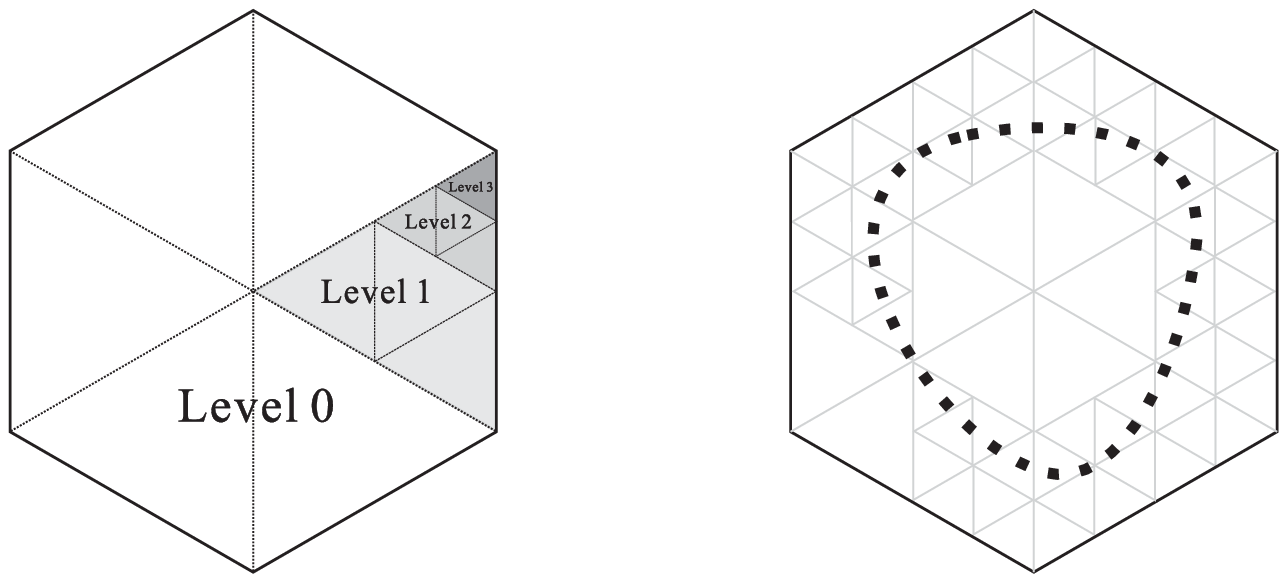}}\\
\small Fig.5 The regular triangle quadtree structure.
\end{figure}

The determination of the relationships between the cells in the triangular tree is quite similar to that in the square tree and so is omitted.

\subsection{Expansions and translations}
We first give the expansions and translations used in FMM and analyze their errors. Then, for square periodic and hexagon periodic arrangements of scatterers, the corresponding periodic tree structure is established, the convergence order of the FMM for different tree structure is given. In the following analysis, the boundary of each scatterer is divided into $N$ elements, the expansions and translations are truncated from $-p$ to $p$.

The multipole expansions, local expansions and translations of the coefficients for $S_1,D_1,K_1,T_1$ have been given in \cite{7}. Thus, we only consider the expansions and translations for $\overline{S}_2,\overline{D}_2,\overline{K}_2$ and $\overline{T}_2$. The following Graf's addition theorem \cite{22} will be used.
\begin{equation}
\mathscr{B}_m(|\mathbf{x}-\mathbf{y}|)e^{\pm{\mathrm{i}m}\theta_{\mathbf{x}-\mathbf{y}}}=\sum_{n=-\infty}^{\infty}\mathscr{B}_{m+n}(|\mathbf{x}|)e^{\pm{\mathrm{i}(m+n)}
\theta_\mathbf{x}}J_n(|\mathbf{y}|)e^{\mp{\mathrm{i}n}\theta_\mathbf{y}},\qquad |\mathbf{y}|<|\mathbf{x}|,
\end{equation}
where $m\in\mathds{Z}$, $\mathscr{B}$ denotes $J,Y,H^{(1)}$ or $H^{(2)}$, $\theta_\mathbf{x-y}$ is the angle between $\mathbf{x-y}$ and the $x$ axis. When $\mathscr{B}=J$, the restriction $|\mathbf{y}|<|\mathbf{x}|$ is unnecessary. We denote the remainder term of (10) as
$$R^{\mathscr{B}}_{m,p}(\mathbf{x},\mathbf{y}):=\left(\sum_{n=p+1}^\infty+\sum_{n=-\infty}^{-p-1}\right)\mathscr{B}_{m+n}(|\mathbf{x}|)
e^{\pm{\mathrm{i}(m+n)}\theta_\mathbf{x}}J_n(|\mathbf{y}|)e^{\mp{\mathrm{i}n}\theta_\mathbf{y}}.$$
For convenience, we let
$$\mathcal{H}^{\pm}_n(\mathbf{x}):=H^{(1)}_n(|\mathbf{x}|)e^{\pm{\mathrm{i}n}\theta_{\mathbf{x}}},\qquad \mathcal{J}^{\pm}_n(\mathbf{x}):=J_n(|\mathbf{x}|)e^{\pm{\mathrm{i}n}\theta_{\mathbf{x}}},\qquad n\in\mathds{Z}.$$

First consider the expansions and translations for $\overline{S}_2$. If the boundary $\Gamma_i$ is divided into $N$ elements $\Delta\Gamma_j,j=1,2,\cdots,N$, then we have
\begin{eqnarray*}
{\overline{S}}_2\varphi_i(\mathbf{x})\hspace{-0.6cm}&&=\int_{\Gamma_i}\sum_{m=-\infty}^\infty\alpha^m\Phi_2(\mathbf{x},\mathbf{y}+m\mathbf{h})\varphi_i(\mathbf{y})ds(\mathbf{y})\\
&&=\frac{\mathrm{i}}{4}\sum_{j=1}^N\int_{\Delta\Gamma_j}\sum_{m=-\infty}^\infty\alpha^mH^{(1)}_0(k|\mathbf{x}-\mathbf{y}-m\mathbf{h}|)\varphi_i(\mathbf{y})ds(\mathbf{y}),
\end{eqnarray*}
where $\mathbf{y}\in\Delta\Gamma_j$. Suppose the tree structure for $\Gamma_i$ has been constructed and $\mathcal{D}$ is a cell of the tree that covers the source point $\mathbf{x}$, the cell $\mathcal{C}$ with centroid $\mathbf{O}_{\mathcal{C}}$ is a well separated cell of $\mathcal{D}$ and $\mathbf{y}\in\Delta\Gamma_j\subset\mathcal{C}$. From (10), when $|\mathbf{y}-\mathbf{O}_\mathcal{C}|<|\mathbf{x}-\mathbf{O}_{\mathcal{C}}-m\mathbf{h}|$, we obtain the following multipole expansion (ME):
\begin{eqnarray}
&&\int_{\Delta\Gamma_j}\sum_{m=-\infty}^\infty\alpha^mH^{(1)}_0(k|\mathbf{x}-\mathbf{y}-m\mathbf{h}|)\varphi_i(\mathbf{y})ds(\mathbf{y})\nonumber\\
&&=\int_{\Delta\Gamma_j}\sum_{m=-\infty}^\infty\alpha^m\left(\sum_{n=-\infty}^{\infty}\mathcal{H}_n^+\big(k(\mathbf{x}-\mathbf{O}_{\mathcal{C}}-m\mathbf{h})\big)
\mathcal{J}^-_n\big(k(\mathbf{y}-\mathbf{O}_{\mathcal{C}})\big)\right)\varphi_i(\mathbf{y})ds(\mathbf{y})\nonumber\\
&&=\sum_{n=-p}^p\sum_{m=-\infty}^\infty\alpha^m\mathcal{H}_n^+\big(k(\mathbf{x}-\mathbf{O}_{\mathcal{C}}-m\mathbf{h})\big)\mathbf{M}_{n}(\mathbf{O}_{\mathcal{C}})+E_{ME}(\mathbf{x},p),
\end{eqnarray}
where
$$\mathbf{M}_n(\mathbf{O}_{\mathcal{C}}):=\int_{\Delta\Gamma_j}\mathcal{J}^-_n\big(k(\mathbf{y}-\mathbf{O}_{\mathcal{C}})\big)\varphi_i(\mathbf{y})ds(\mathbf{y}),
\qquad\mathbf{y}\in\Delta\Gamma_j$$
is the multipole moment (MM) and
\begin{equation}
E_{ME}(\mathbf{x},p):=\int_{\Delta\Gamma_j}\sum_{m=-\infty}^\infty\alpha^mR_{0,p}^{H}\big(k(\mathbf{x}-\mathbf{O}_{\mathcal{C}}-m\mathbf{h}),k(\mathbf{y}-\mathbf{O}_{\mathcal{C}})\big)
\varphi_i(\mathbf{y})ds(\mathbf{y}).
\end{equation}

Let $\mathcal{CS}$ denotes a child of $\mathcal{C}$ and $\mathbf{O}_\mathcal{CS}$ is the centroid of $\mathcal{CS}$. From (10), for each $\mathbf{y}\in\Delta\Gamma_j\subset\mathcal{CS}$, we have the following M2M translation:
\begin{eqnarray*}
\mathbf{M}_n(\mathbf{O}_{\mathcal{C}})\hspace{-0.6cm}&&=\int_{\Delta\Gamma_j}\mathcal{J}^-_n\big(k(\mathbf{y}-\mathbf{O}_{\mathcal{C}})\big)\varphi_i(\mathbf{y})ds(\mathbf{y})\\
&&=\int_{\Delta\Gamma_j}\left(\sum_{l=-\infty}^{\infty}\mathcal{J}^-_l\big(k(\mathbf{y}-\mathbf{O}_\mathcal{CS})\big)\mathcal{J}^-_{n-l}\big(k(\mathbf{O}_\mathcal{CS}
-\mathbf{O}_\mathcal{C})\big)\right)\varphi_i(\mathbf{y})ds(\mathbf{y})\\
&&=\sum_{l=-p}^{p}\mathbf{M}_l(\mathbf{O}_\mathcal{CS})\mathcal{J}^-_{n-l}\big(k(\mathbf{O}_\mathcal{CS}-\mathbf{O}_\mathcal{C})\big)+E_{M2M}(\mathbf{O}_\mathcal{C},n,p),
\end{eqnarray*}
where
$$E_{M2M}(\mathbf{O}_\mathcal{C},n,p):=\int_{\Delta\Gamma_j}R_{n,p}^{J}\big(k(\mathbf{O}_\mathcal{CS}-\mathbf{O}_\mathcal{C}),
k(\mathbf{O}_\mathcal{CS}-\mathbf{y})\big)\varphi_i(\mathbf{y})ds(\mathbf{y}).$$

For the main part of the multipole expansion (11), when $|\mathbf{x}-\mathbf{O}_\mathcal{D}|<|\mathbf{O}_\mathcal{D}-\mathbf{O}_\mathcal{C}-m\mathbf{h}|$,
we obtain the local expansion (LE):
\begin{eqnarray*}
&&\hspace{-0.6cm}\sum_{n=-p}^p\sum_{m=-\infty}^\infty\alpha^m\mathcal{H}_n^+\big(k(\mathbf{x}-\mathbf{O}_{\mathcal{C}}-m\mathbf{h})\big)\mathbf{M}_{n}(\mathbf{O}_{\mathcal{C}})\\
&&\hspace{-0.6cm}=\sum_{n=-p}^p\sum_{m=-\infty}^\infty\alpha^m\left(\sum_{l=-\infty}^{\infty}\mathcal{J}^+_l\big(k(\mathbf{x}-\mathbf{O}_\mathcal{D})\big)
\mathcal{H}^+_{n-l}\big(k(\mathbf{O}_\mathcal{D}-\mathbf{O}_\mathcal{C}-m\mathbf{h})\big)\right)\mathbf{M}_{n}(\mathbf{O}_{\mathcal{C}})\\
&&\hspace{-0.6cm}=\sum_{l=-p}^p\mathbf{L}_l(\mathbf{O}_\mathcal{D})\mathcal{J}^+_l\big(k(\mathbf{x}-\mathbf{O}_\mathcal{D})\big)+E_{M2L}(\mathbf{x},\mathbf{O}_\mathcal{D},p),
\end{eqnarray*}
where
$$\mathbf{L}_l(\mathbf{O}_\mathcal{D}):=\sum_{n=-p}^p\sum_{m=-\infty}^\infty\alpha^m\mathcal{H}^+_{n-l}\big(k(\mathbf{O}_\mathcal{D}-\mathbf{O}_\mathcal{C}-m\mathbf{h})\big)
\mathbf{M}_{n}(\mathbf{O}_{\mathcal{C}})$$
is the M2L translation and
\begin{equation}
E_{M2L}(\mathbf{x},\mathbf{O}_\mathcal{D},p):=\sum_{n=-p}^p\sum_{m=-\infty}^\infty\alpha^mR_{n,p}^{H}\big(k(\mathbf{O}_\mathcal{D}-\mathbf{O}_\mathcal{C}
-m\mathbf{h}),k(\mathbf{O}_\mathcal{D}-\mathbf{x})\big)\mathbf{M}_{n}(\mathbf{O}_{\mathcal{C}}).
\end{equation}

Suppose $\mathcal{DP}$ is the parent of $\mathcal{D}$ and $\mathbf{O}_\mathcal{DP}$ is the centroid of $\mathcal{DP}$. From (10), for each $\mathbf{x}\in\mathcal{D}\subset\mathcal{DP}$, we have
\begin{eqnarray*}
&&\sum_{l=-p}^p\mathbf{L}_l(\mathbf{O}_\mathcal{DP})\mathcal{J}^+_l\big(k(\mathbf{x}-\mathbf{O}_\mathcal{DP})\big)\\
&&=\sum_{l=-p}^{p}\mathbf{L}_l(\mathbf{O}_\mathcal{DP})\left(\sum_{n=-\infty}^{\infty}\mathcal{J}^+_n\big(k(\mathbf{x}-\mathbf{O}_\mathcal{D})\big)\mathcal{J}^+_{l-n}
\big(k(\mathbf{O}_\mathcal{D}-\mathbf{O}_\mathcal{DP})\big)\right)\\
&&=\sum_{n=-p}^{p}\mathbf{L}_n(\mathbf{O}_\mathcal{D})\mathcal{J}^+_n\big(k(\mathbf{x}-\mathbf{O}_\mathcal{D})\big)+E_{L2L}(\mathbf{x},\mathbf{O}_\mathcal{D},p),
\end{eqnarray*}
where
$$\mathbf{L}_n(\mathbf{O}_\mathcal{D})=\sum_{l=-p}^{p}\mathbf{L}_l(\mathbf{O}_\mathcal{DP})\mathcal{J}^+_{l-n}\big(k(\mathbf{O}_\mathcal{D}-\mathbf{O}_\mathcal{DP})\big)$$
is the L2L translation and
$$E_{L2L}(\mathbf{x},\mathbf{O}_\mathcal{D},p):=\sum_{l=-p}^{p}\mathbf{L}_l(\mathbf{O}_\mathcal{DP})R_{l,p}^{J}\big(k(\mathbf{O}_\mathcal{D}-\mathbf{O}_\mathcal{DP}),
k(\mathbf{O}_\mathcal{D}-\mathbf{x})\big).$$

In the above analysis, $E_{ME},E_{M2M},E_{M2L},E_{L2L}$ denote the truncation errors of the ME, M2M, M2L and L2L respectively.

From \cite{23}, we see that $E_{M2M}$ and $E_{L2L}$ is more smaller than $E_{ME}$ and $E_{M2L}$, thus we only give the estimation of the convergence of $E_{ME}$ and $E_{M2L}$.

{\theorem Let $\mathcal{D}$ be a cell of the tree structure and let $\mathcal{C}$ be a well separated cell of $\mathcal{D}$. For each source point $\mathbf{x}\in\mathcal{D}$, when
$kH(1\pm\sin\theta)\neq2n\pi$ and $H\geq|\mathbf{x}-\mathbf{O}_{\mathcal{C}}|$,
$$\big|E_{ME}(\mathbf{x},p)\big|=\mathcal{O}\left(\frac{\gamma^p}{p}\right),$$
where}
$$\gamma:=\max_{\mathbf{y}\in\Delta\Gamma_j}\left\{\left.\frac{|\mathbf{y}-\mathbf{O}_{\mathcal{C}}|}{|\mathbf{x}-\mathbf{O}_{\mathcal{C}}-m\mathbf{h}|}\right|m=-1,0,1\right\}.$$
\textbf{Proof.} It can be seen from (12),
\begin{equation}
\big|E_{ME}(\mathbf{x},p)\big|\leq\int_{\Delta\Gamma_j}\sum_{m=-\infty}^\infty\left|R_{0,p}^{H}\big(k(\mathbf{x}-\mathbf{O}_{\mathcal{C}}-m\mathbf{h}),k(\mathbf{y}
-\mathbf{O}_{\mathcal{C}})\big)\right||\varphi_i(\mathbf{y})|ds(\mathbf{y}),
\end{equation}
In \cite{23}, we have proved that
\begin{equation}
\left|R_{0,p}^{H}\big(k(\mathbf{x}-\mathbf{O}_{\mathcal{C}}-m\mathbf{h}),k(\mathbf{y}-\mathbf{O}_{\mathcal{C}})\big)\right|=\mathcal{O}\left(\frac{\gamma_m^p}{p}\right),
\end{equation}
where
$$\gamma_m:=\frac{|\mathbf{y}-\mathbf{O}_{\mathcal{C}}|}{|\mathbf{x}-\mathbf{O}_{\mathcal{C}}-m\mathbf{h}|}.$$

When $m\geq2$ and $H\geq|\mathbf{x}-\mathbf{O}_{\mathcal{C}}|$,
$$\gamma_{\pm{m}}=\frac{|\mathbf{y}-\mathbf{O}_{\mathcal{C}}|}{|\mathbf{x}-\mathbf{O}_{\mathcal{C}}\pm{m}\mathbf{h}|}\leq\frac{|\mathbf{y}-\mathbf{O}_{\mathcal{C}}|}
{mH-|\mathbf{x}-\mathbf{O}_{\mathcal{C}}|}\leq\frac{|\mathbf{y}-\mathbf{O}_{\mathcal{C}}|}{(m-1)|\mathbf{x}-\mathbf{O}_{\mathcal{C}}|}=\frac{\gamma_0}{m-1},$$
thus, we have
$$\sum_{m=-\infty}^\infty\gamma_m^p\leq\gamma_{-1}^p+\gamma_1^p+\left(1+2\sum_{m=2}^\infty\frac{1}{(m-1)^p}\right)\gamma_0^p,$$
since when $p\geq2$,
$$\sum_{m=2}^\infty\frac{1}{(m-1)^p}\leq\frac{\pi^2}{6},$$
it follows that
\begin{equation}
\sum_{m=-\infty}^\infty\frac{\gamma_m^p}{p}=\mathcal{O}\left(\frac{\gamma^p}{p}\right),
\end{equation}
where $\gamma=\max\{\gamma_{-1},\gamma_0,\gamma_1\}$. From (14) (15) and (16), we prove the theorem. $~~~\Box$

The same proof remains valid for $E_{M2L}(\mathbf{x},\mathbf{O}_\mathcal{D},p)$. From (13),
$$\big|E_{M2L}(\mathbf{x},\mathbf{O}_\mathcal{D},p)\big|\leq\sum_{m=-\infty}^\infty\left|\sum_{n=-p}^pR_{n,p}^{H}\big(k(\mathbf{O}_\mathcal{D}-\mathbf{O}_\mathcal{C}
-m\mathbf{h}),k(\mathbf{O}_\mathcal{D}-\mathbf{x})\big)\mathbf{M}_{n}(\mathbf{O}_{\mathcal{C}})\right|,$$
by \cite{23}, we have
$$\left|\sum_{n=-p}^pR_{n,p}^{H}\big(k(\mathbf{O}_\mathcal{D}-\mathbf{O}_\mathcal{C}-m\mathbf{h}),k(\mathbf{O}_\mathcal{D}-\mathbf{x})\big)\mathbf{M}_{n}
(\mathbf{O}_{\mathcal{C}})\right|=\mathcal{O}\left(\frac{\lambda_m^p}{p}\right),$$
where
$$\lambda_m:=\frac{|\mathbf{x}-\mathbf{O}_\mathcal{D}|}{|\mathbf{O}_\mathcal{D}-\mathbf{O}_\mathcal{C}-m\mathbf{h}|}\exp\left({\frac{3|\mathbf{y}-\mathbf{O}_{\mathcal{C}}|}
{2|\mathbf{O}_\mathcal{D}-\mathbf{O}_\mathcal{C}-m\mathbf{h}|}}\right).$$
We can also derive the estimation of the convergence order of $E_{M2L}(\mathbf{x},\mathbf{O}_\mathcal{D},p)$ as follows.

{\theorem For each source point $\mathbf{x}\in\mathcal{D}$, when $kH(1\pm\sin\theta)\neq2n\pi$ and $H\geq|\mathbf{O}_\mathcal{D}-\mathbf{O}_\mathcal{C}|$,
$$\big|E_{M2L}(\mathbf{x},\mathbf{O}_\mathcal{D},p)\big|=\mathcal{O}\left(\frac{\lambda^p}{p}\right),$$
where $\lambda:=\max\left\{\lambda_{-1},\lambda_0,\lambda_1\right\}$.}\\
\textbf{Proof.} When $m\geq2$ and $H\geq|\mathbf{O}_\mathcal{D}-\mathbf{O}_\mathcal{C}|$,
$$\lambda_{\pm{m}}\leq\frac{|\mathbf{x}-\mathbf{O}_\mathcal{D}|}{(m-1)|\mathbf{O}_\mathcal{D}-\mathbf{O}_\mathcal{C}|}\exp\left({\frac{3|\mathbf{y}-\mathbf{O}_{\mathcal{C}}|}
{2|\mathbf{O}_\mathcal{D}-\mathbf{O}_\mathcal{C}|}}\right)=\frac{\lambda_0}{m-1},$$
which proves the theorem. $~~~\Box$

The above two theorems show that the convergence order of the error of FMM for $\overline{S}_2$ is $\mathcal{O}(\tau^p/p)$, where $\tau=\max\left\{\gamma,\lambda\right\}$.

For the square and triangular quad-tree structures proposed in Section 3.1, the values of $\gamma$ and $\lambda$ can be easily calculated and are shown in Table 1, where $\mathcal{L}=L_\mathcal{D}-L_\mathcal{C}$, $L_\mathcal{D}$ and $L_\mathcal{C}$ are the layer numbers of cells $\mathcal{D}$ and $\mathcal{C}$ respectively, $\mathcal{L}=0$ corresponds to the symmetric tree and $\mathcal{L}>0$ corresponds to the asymmetric tree.
\begin{center}
{Table 1. Values of $\gamma, \lambda, \tau$ for two tree structures.}\\
{\begin{tabular}{p{1.2cm}<{\centering}|p{1.5cm}<{\centering}|p{1.5cm}<{\centering}|p{1.5cm}<{\centering}|p{1.5cm}<{\centering}|p{1.5cm}<{\centering}}\hline
& \multicolumn{3}{c|}{square tree} & \multicolumn{2}{c}{triangular tree} \\\hline
& $\mathcal{L}=0$ & $\mathcal{L}=1$ & $\mathcal{L}=2$ & $\mathcal{L}=0$ & $\mathcal{L}=1$\\\hline
 $\gamma$ & $0.4714$ & $0.7071$ & $0.9428$ & $0.4000$ & $0.7559$ \\\hline
 $\lambda$ & $0.6009$ & $0.6374$ & $0.6640$ & $0.6663$ & $0.6863$ \\\hline
 $\tau$ & $0.6009$ & $0.7071$ & $0.9428$ & $0.6663$ & $0.7559$ \\\hline
\end{tabular}}
\end{center}

The derivations of the expansions and translations for $\overline{D}_2,\overline{K}_2,\overline{T}_2$ are quite similar to that for $\overline{S}_2$ and so is omitted. Many formulas for $\overline{D}_2,\overline{K}_2,\overline{T}_2$ are the same as that for $\overline{S}_2$, thus we only show the different parts as follows.

\textbf{Expansions and moments for} $\overline{D}_2$:
\begin{eqnarray*}
&&\mathrm{MM}:~~\mathbf{M}_n(\mathbf{O}_\mathcal{C})=\int_{\Delta\Gamma_j}\frac{\partial\mathcal{J}^-_n\big(k(\mathbf{y}-\mathbf{O}_\mathcal{C})\big)}{\partial{\nu(\mathbf{y})}}
\phi_i(\mathbf{y})ds(\mathbf{y}),\\
&&E_{ME}(\mathbf{x},p)=\int_{\Delta\Gamma_j}\sum_{m=-\infty}^\infty\alpha^m\frac{\partial{R}_{0,p}^{H}\big(k(\mathbf{x}-\mathbf{O}_{\mathcal{C}}-m\mathbf{h}),k(\mathbf{y}-\mathbf{O}_{\mathcal{C}})
\big)}{\partial{\nu(\mathbf{y})}}\varphi_i(\mathbf{y})ds(\mathbf{y}),\hspace{0.9cm}\\
&&E_{M2M}(\mathbf{O}_\mathcal{C},n,p)=\int_{\Delta\Gamma_j}\frac{\partial{R}_{n,p}^{J}\big(k(\mathbf{O}_\mathcal{CS}-\mathbf{O}_\mathcal{C}),
k(\mathbf{O}_\mathcal{CS}-\mathbf{y})\big)}{\partial{\nu(\mathbf{y})}}\varphi_i(\mathbf{y})ds(\mathbf{y}).
\end{eqnarray*}

\textbf{Expansions and moments for} $\overline{K}_2$:
\begin{eqnarray*}
&&\mathrm{ME}:~~\sum_{n=-p}^{p}\sum_{m=-\infty}^\infty\alpha^m\frac{\partial\mathcal{H}^+_n\big(k(\mathbf{x}-\mathbf{O}_\mathcal{C}-m\mathbf{h})\big)}
{\partial{\nu(\mathbf{x})}}{\mathbf{M}}_{n}(\mathbf{O}_\mathcal{C}),\\
&&\mathrm{MM}:~~\mathbf{M}_n(\mathbf{O}_{\mathcal{C}})=\int_{\Delta\Gamma_j}\mathcal{J}^-_n\big(k(\mathbf{y}-\mathbf{O}_\mathcal{C})\big)\varphi_i(\mathbf{y})ds(\mathbf{y}),\\
&&\mathrm{LE}:~~\sum_{l=-p}^{p}\mathbf{L}_l(\mathbf{O}_\mathcal{D})\frac{\partial\mathcal{J}^+_l\big(k(\mathbf{x}-\mathbf{O}_\mathcal{D})\big)}{\partial{\nu(\mathbf{x})}},\\
&&E_{ME}(\mathbf{x},p)=\int_{\Delta\Gamma_j}\sum_{m=-\infty}^\infty\alpha^m\frac{\partial{R}_{0,p}^{H}\big(k(\mathbf{x}-\mathbf{O}_{\mathcal{C}}-m\mathbf{h}),k(\mathbf{y}-\mathbf{O}_{\mathcal{C}})
\big)}{\partial{\nu(\mathbf{x})}}\varphi_i(\mathbf{y})ds(\mathbf{y}),\\
&&E_{M2L}(\mathbf{x},\mathbf{O}_\mathcal{D},p)=\sum_{n=-p}^p\sum_{m=-\infty}^\infty\alpha^m\frac{\partial{R}_{n,p}^{H}\big(k(\mathbf{O}_\mathcal{D}-\mathbf{O}_\mathcal{C}
-m\mathbf{h}),k(\mathbf{O}_\mathcal{D}-\mathbf{x})\big)}{\partial{\nu(\mathbf{x})}}\mathbf{M}_{n}(\mathbf{O}_{\mathcal{C}}),\\
&&E_{L2L}(\mathbf{x},\mathbf{O}_\mathcal{D},p)=\sum_{l=-p}^{p}\mathbf{L}_l(\mathbf{O}_\mathcal{DP})\frac{\partial{R}_{l,p}^{J}\big(k(\mathbf{O}_\mathcal{D}-\mathbf{O}
_\mathcal{DP}),k(\mathbf{O}_\mathcal{D}-\mathbf{x})\big)}{\partial{\nu(\mathbf{x})}}.
\end{eqnarray*}

\textbf{Expansions and moments for} $\overline{T}_2$:
\begin{eqnarray*}
&&\mathrm{ME}:~~\sum_{n=-p}^{p}\sum_{m=-\infty}^\infty\alpha^m\frac{\partial\mathcal{H}^+_n\big(k(\mathbf{x}-\mathbf{O}_\mathcal{C}-m\mathbf{h})\big)}
{\partial{\nu(\mathbf{x})}}{\mathbf{M}}_{n}(\mathbf{O}_\mathcal{C}),\\
&&\mathrm{MM}:~~\mathbf{M}_n(\mathbf{O}_\mathcal{C})=\int_{\Delta\Gamma_j}\frac{\partial\mathcal{J}^-_n\big(k(\mathbf{y}-\mathbf{O}_\mathcal{C})\big)}
{\partial{\nu(\mathbf{y})}}\phi_i(\mathbf{y})ds(\mathbf{y}),\\
&&\mathrm{LE}:~~\sum_{l=-p}^{p}\mathbf{L}_l(\mathbf{O}_\mathcal{D})\frac{\partial\mathcal{J}^+_l\big(k(\mathbf{x}-\mathbf{O}_\mathcal{D})\big)}{\partial{\nu(\mathbf{x})}},\\
&&E_{ME}(\mathbf{x},p)=\frac{\partial}{\partial{\nu(\mathbf{x})}}\int_{\Delta\Gamma_j}\sum_{m=-\infty}^\infty\alpha^m\frac{\partial{R}_{0,p}^{H}\big(k(\mathbf{x}
-\mathbf{O}_{\mathcal{C}}-m\mathbf{h}),k(\mathbf{y}-\mathbf{O}_{\mathcal{C}})
\big)}{\partial{\nu(\mathbf{y})}}\varphi_i(\mathbf{y})ds(\mathbf{y}),\\
&&E_{M2M}(\mathbf{O}_\mathcal{C},n,p)=\int_{\Delta\Gamma_j}\frac{\partial{R}_{n,p}^{J}\big(k(\mathbf{O}_\mathcal{CS}-\mathbf{O}_\mathcal{C}),
k(\mathbf{O}_\mathcal{CS}-\mathbf{y})\big)}{\partial{\nu(\mathbf{y})}}\varphi_i(\mathbf{y})ds(\mathbf{y}),\\
&&E_{M2L}(\mathbf{x},\mathbf{O}_\mathcal{D},p)=\sum_{n=-p}^p\sum_{m=-\infty}^\infty\alpha^m\frac{\partial{R}_{n,p}^{H}\big(k(\mathbf{O}_\mathcal{D}-\mathbf{O}_\mathcal{C}
-m\mathbf{h}),k(\mathbf{O}_\mathcal{D}-\mathbf{x})\big)}{\partial{\nu(\mathbf{x})}}\mathbf{M}_{n}(\mathbf{O}_{\mathcal{C}}),\\
&&E_{L2L}(\mathbf{x},\mathbf{O}_\mathcal{D},p)=\sum_{l=-p}^{p}\mathbf{L}_l(\mathbf{O}_\mathcal{DP})\frac{\partial{R}_{l,p}^{J}\big(k(\mathbf{O}_\mathcal{D}-\mathbf{O}
_\mathcal{DP}),k(\mathbf{O}_\mathcal{D}-\mathbf{x})\big)}{\partial{\nu(\mathbf{x})}}.
\end{eqnarray*}

In addition, from the recurrence relation \cite{22}
$$2\mathscr{B}'_n(z)=\mathscr{B}_{n-1}(z)-\mathscr{B}_{n+1}(z),$$
we can also derive that the convergence orders of the errors of FMM for $\overline{D}_2,\overline{K}_2,\overline{T}_2$ are $\mathcal{O}(\tau^p)$, $\mathcal{O}(\tau^p)$ and $\mathcal{O}(p\tau^p)$ respectively.

\subsection{Fast computation of the M2L translation}
Lemma 1 shows that a larger truncation number $q$ is needed to ensure the accuracy of the algorithm. This leads to a huge amount of computation for the M2L translation:
$$\mathbf{L}_l(\mathbf{O}_\mathcal{D})=\sum_{n=-p}^p\left(\sum_{m=-\infty}^\infty\alpha^mH_{l-n}^{(1)}(k|\mathbf{O}_\mathcal{D}-\mathbf{O}_\mathcal{C}-m\mathbf{h}|)
e^{\mathrm{i}(l-n)\theta_{\mathbf{O}_\mathcal{D}-\mathbf{O}_\mathcal{C}-m\mathbf{h}}}\right)\mathbf{M}_{n}(\mathbf{O}_{\mathcal{C}}).$$
We will show the fast computation of the sum:
$$\sum_{m=-\infty}^\infty\alpha^mH_n^{(1)}(k|\mathbf{O}_\mathcal{D}-\mathbf{O}_\mathcal{C}-m\mathbf{h}|)e^{\mathrm{i}n\theta_{\mathbf{O}_\mathcal{D}-\mathbf{O}_\mathcal{C}
-m\mathbf{h}}}.\qquad n\in\mathds{Z}.$$

In \cite{18}, the Fourier transform was applied to accelerate the computation of the periodic Green's function for Helmholtz equation. For each integer $q\geq1$ and $\mathbf{x},\mathbf{y}\in\mathds{R}^2$,
$$\sum_{m=q+1}^{+\infty}e^{\mathrm{i}m\beta}H_0^{(1)}(k|\mathbf{x}-\mathbf{y}-m\mathbf{h}|)=\frac{2}{\mathrm{i}}\int_{-\infty}^{+\infty}\frac{e^{\mathrm{i}q\beta+\mathrm{i}(x_1-y_1)t
+(x_2-y_2-qH)\mu}}{\mu(e^{H\mu-\mathrm{i}\beta}-1)}dt,$$
$$\sum_{m=-\infty}^{-q-1}e^{\mathrm{i}m\beta}H_0^{(1)}(k|\mathbf{x}-\mathbf{y}-m\mathbf{h}|)=\frac{2}{\mathrm{i}}\int_{-\infty}^{+\infty}\frac{e^{-\mathrm{i}q\beta+\mathrm{i}(x_1-y_1)t
+(y_2-x_2-qH)\mu}}{\mu(e^{H\mu+\mathrm{i}\beta}-1)}dt,$$
where $\beta=kH\sin\theta$, $\mathbf{x}=(x_1,x_2),\mathbf{y}=(y_1,y_2),\mathbf{h}=(0,H)$ and $\mu=\sqrt{t^2-k^2}$. Since for each $n\geq0$,
$$H_n^{(1)}(k|\mathbf{x}|)e^{\mathrm{i}n\theta_\mathbf{x}}=\frac{1}{(-k)^n}\left(\frac{\partial}{\partial{x_1}}+\mathrm{i}\frac{\partial}{\partial{x_2}}\right)^n
H_0^{(1)}(k|\mathbf{x}|),$$
$$H_{-n}^{(1)}(k|\mathbf{x}|)e^{-\mathrm{i}n\theta_\mathbf{x}}=\frac{1}{k^n}\left(\frac{\partial}{\partial{x_1}}-\mathrm{i}\frac{\partial}{\partial{x_2}}\right)^n
H_0^{(1)}(k|\mathbf{x}|),$$
it follows that
$$\sum_{m=q+1}^{+\infty}e^{\mathrm{i}m\beta}{H_{\pm{n}}^{(1)}}(k|\mathbf{x}-\mathbf{y}-m\mathbf{h}|)e^{\pm\mathrm{i}n\theta_{\mathbf{x}-\mathbf{y}-m\mathbf{h}}}=\frac{2}{\mathrm{i}^{1+n}k^n}
\int_{-\infty}^{+\infty}\frac{e^{\mathrm{i}q\beta+\mathrm{i}(x_1-y_1)t+(x_2-y_2-qH)\mu}(\mu\pm{t})^n}{\mu(e^{H\mu-\mathrm{i}\beta}-1)}dt,$$
$$\sum_{m=-\infty}^{-q-1}e^{\mathrm{i}m\beta}{H_{\pm{n}}^{(1)}}(k|\mathbf{x}-\mathbf{y}-m\mathbf{h}|)e^{\pm\mathrm{i}n\theta_{\mathbf{x}-\mathbf{y}-m\mathbf{h}}}=\frac{2\mathrm{i}^{1+n}}{k^n}
\int_{-\infty}^{+\infty}\frac{e^{-\mathrm{i}q\beta+\mathrm{i}(y_1-x_1)t+(y_2-x_2-qH)\mu}(\mu\pm{t})^n}{\mu(e^{H\mu+\mathrm{i}\beta}-1)}dt.$$
Thus, we have
\begin{eqnarray}
&&\sum_{m=-\infty}^{+\infty}\alpha^m{H_{\pm{n}}^{(1)}}(k|\mathbf{x}-\mathbf{y}-m\mathbf{h}|)e^{\pm\mathrm{i}n\theta_{\mathbf{x}-\mathbf{y}-m\mathbf{h}}}\nonumber\\
=\hspace{-0.6cm}&&\sum_{m=-q}^q\alpha^m{H_{\pm{n}}^{(1)}}(k|\mathbf{x}-\mathbf{y}-m\mathbf{h}|)e^{\pm\mathrm{i}n\theta_{\mathbf{x}-\mathbf{y}-m\mathbf{h}}}\nonumber\\
&&+~\frac{2}{\mathrm{i}^{1+n}k^n}\int_{-\infty}^{+\infty}\frac{e^{\mathrm{i}q\beta+\mathrm{i}(x_1-y_1)t+(x_2-y_2-qH)\mu}(\mu\pm{t})^n}{\mu(e^{H\mu-\mathrm{i}\beta}-1)}dt\nonumber\\
&&+~\frac{2\mathrm{i}^{1+n}}{k^n}\int_{-\infty}^{+\infty}\frac{e^{-\mathrm{i}q\beta+\mathrm{i}(y_1-x_1)t+(y_2-x_2-qH)\mu}(\mu\pm{t})^n}{\mu(e^{H\mu+\mathrm{i}\beta}-1)}dt.
\end{eqnarray}
Since the integrands in (17) has poles on the real axis ($t=\pm{k}$), we use the curve
$$C: t=\pm\sqrt{u^2-2\mathrm{i}ku} \qquad u>0 $$
as the path of integration, hence
\begin{eqnarray}
&&\int_0^{+\infty}\frac{e^{\mathrm{i}q\beta+\mathrm{i}(x_1-y_1)t+(x_2-y_2-qH)\mu}(\mu+t)^n}{\mu(e^{H\mu-\mathrm{i}\beta}-1)}dt\nonumber\\&&=\int_0^{+\infty}\frac{e^{\mathrm{i}q\beta
+\mathrm{i}(x_1-y_1)\sqrt{u^2-2\mathrm{i}ku}+(x_2-y_2-qH)(u-\mathrm{i}k)}\left(u-\mathrm{i}k+\sqrt{u^2-2\mathrm{i}ku}\right)^n}{(u-\mathrm{i}k)(e^{Hu-\mathrm{i}Hk-\mathrm{i}\beta}-1)}du.
\end{eqnarray}
Note that the integrand in (18) is analytic when $kH(1+\sin\theta)\neq2n\pi$.

In this paper, the integrals are computed by the quadrature rule given in \cite{24}, that is
\begin{equation}
\int_0^{+\infty}f(x,u)du\approx\sum_{i=1}^P\omega_if(x,u_i),
\end{equation}
where $\omega_i$ and $u_i$ are weights and nodes respectively. The remainder term of the quadrature rule for (18) can be written as
$$\mathscr{R}_{n,\omega}:=\int_\omega^{+\infty}\frac{e^{\mathrm{i}q\beta+\mathrm{i}(x_1-y_1)\sqrt{u^2-2\mathrm{i}ku}
+(x_2-y_2-qH)(u-\mathrm{i}k)}\left(u-\mathrm{i}k+\sqrt{u^2-2\mathrm{i}ku}\right)^n}{(u-\mathrm{i}k)(e^{Hu-\mathrm{i}Hk-\mathrm{i}\beta}-1)}du,$$
where
$$\omega=\sum_{i=1}^P\omega_i.$$
We give the estimate for the bound on $\mathscr{R}_{n,\omega}$ in the following theorem.

{\theorem Let $kH(1+\sin\theta)\neq2n\pi$ and $\omega\geq\max\{\ln2/H,k\}$, when $q\geq2(n-1)/H\omega$,}
$$|\mathscr{R}_{n,\omega}|\leq\frac{4e^{kL}\left(\sqrt{2}+\sqrt{3}\right)^n}{qH}w^{n-1}e^{-qH\omega}.$$
\textbf{Proof.} Let $\sqrt{u^2-2\mathrm{i}ku}=A+\mathrm{i}B$, an easy computation shows that
$$A^2=\frac{u^2+\sqrt{u^4+4k^2u^2}}{2},\qquad B^2=\frac{2k^2u^2}{u^2+\sqrt{u^4+4k^2u^2}},$$
thus we can derive that $|A|\leq\sqrt{u^2+k^2}$ and $|B|\leq{k}$.

Since $|x_2-y_2|<H$ and $|x_1-y_1|<L$, it follows that
\begin{eqnarray*}
&&\left|\frac{e^{\mathrm{i}q\beta+\mathrm{i}(x_1-y_1)\sqrt{u^2-2\mathrm{i}ku}+(x_2-y_2-qH)(u-\mathrm{i}k)}\left(u-\mathrm{i}k+\sqrt{u^2-2\mathrm{i}ku}\right)^n}{(u-\mathrm{i}k)
(e^{Hu-\mathrm{i}Hk-\mathrm{i}\beta}-1)}\right|\\
&&\leq\frac{e^{B(y_1-x_1)+(x_2-y_2)u}\left(\sqrt{u^2+k^2}+\sqrt{A^2+B^2}\right)^n}{e^{qHu}(e^{Hu}-1)\sqrt{u^2+k^2}}\\
&&\leq\frac{2e^{kL}\left(\sqrt{u^2+k^2}+\sqrt{u^2+2k^2}\right)^n}{e^{qHu}\sqrt{u^2+k^2}}\\
&&\leq\frac{2\left(\sqrt{2}+\sqrt{3}\right)^ne^{kL}u^{n-1}}{e^{qHu}}
\end{eqnarray*}
In the above estimate, $Hu\geq\ln2$ and $u\geq{k}$ is used. Thus, when $\omega\geq\max\{\ln2/H,k\}$, we have
\begin{equation}
|\mathscr{R}_{n,\omega}|\leq2\left(\sqrt{2}+\sqrt{3}\right)^ne^{kL}\int_\omega^{+\infty}u^{n-1}e^{-qHu}du=2e^{kL}\left(\frac{\sqrt{2}+\sqrt{3}}{qH}\right)^n\Gamma(n,qH\omega),
\end{equation}
where $\Gamma(n,\cdot)$ is the incomplete gamma function. From \cite{22}, we see that
$$\Gamma(n,qH\omega)=(qH\omega)^{n-1}e^{-qH\omega}\left(\sum_{j=0}^{m-1}\frac{(n-1)(n-2)\cdots(n-j)}{(qH\omega)^j}+\varepsilon_m(n,qH\omega)\right),$$
and when $m\geq{n}-1$,
$$|\varepsilon_m(n,qH\omega)|\leq\frac{|(n-1)(n-2)\cdots(n-m)|}{(qH\omega)^m}.$$
Setting $m=n-1$, we derive
\begin{equation}
\left|\Gamma(n,qH\omega)\right|\leq(qH\omega)^{n-1}e^{-qH\omega}\sum_{j=0}^{n-1}\frac{(n-1)(n-2)\cdots(n-j)}{(qH\omega)^j},
\end{equation}
and when $qH\omega\geq2(n-1)$,
\begin{equation}
\sum_{j=0}^{n-1}\frac{(n-1)(n-2)\cdots(n-j)}{(qH\omega)^j}\leq\sum_{j=0}^{n-1}\left(\frac{n-1}{qH\omega}\right)^j<\frac{qH\omega}{qH\omega-n+1}\leq2.
\end{equation}
Form (20), (21) and (22), we prove the theorem. $~~~\Box$

It is worth mentioning that, with the exception of $\mathscr{R}_{n,\omega}$, the error of quadrature rule (19) also includes the following term:
$$\int_0^\omega{f}(x,u)du-\sum_{i=1}^P\omega_if(x,u_i).$$
However, this error depends only on the selection of $\omega_i$ and $u_i$, and which can provide high accuracy.

Theorem 3 shows that the treatment (17) is more efficient than the direct calculation $\mathcal{O}(1/\sqrt{q})$. However, when $n$ is large, the error decays slowly, it follows that a larger $q$ is needed to meet the the accuracy requirement. For a given tolerance error $\epsilon$, we can determine $q$ by the following formula:
\begin{equation}
q=\max\left\{\left[\frac{\ln4+kL+n\ln(\sqrt{2}+\sqrt{3})\omega-\ln\epsilon}{H\omega}\right]+1,\left[\frac{2n-2}{H\omega}\right]+1\right\}.
\end{equation}

\section{Applications in liquid phononic crystals}
It is well known that the acoustic wave propagation in liquid phononic crystal satisfies the boundary value problem (1)-(3). In this section, the proposed periodic FM-BEM will be used to compute the total acoustic field $u(\mathbf{x})$ in the water ($\rho=1000\mathrm{kg}/\mathrm{m^3},c=1500\mathrm{m}/\mathrm{s}$) and mercury ($\rho=13600\mathrm{kg}/\mathrm{m^3},c=1450\mathrm{m}/\mathrm{s}$) phononic crystals. We use the preconditioned GMRES to solve the discretized integral equations, and the numerical experiments were conducted using a FORTRAN program, where the truncation number $p=10$ and $q$ is calculated by (23).

Suppose the incident wave $u^{inc}(\mathbf{x})$ is incident horizontally from the left side, i.e. $\mathbf{d}=(1,0)$. We first compute the scattered field $u^s(\mathbf{x})$ on the right side ($x_1=L$) and then derive the energy transmission coefficient as follows:
$$ETC=\frac{1}{kH}\int_0^H\left|\frac{\partial{u}}{\partial{x}_1}(L,x_2)\right|dx_2,$$
where $\mathbf{x}=(x_1,x_2)$ and $u(\mathbf{x})=u^{inc}(\mathbf{x})+u^s(\mathbf{x})$. We calculate $ETC$ for different incident frequencies, when $ETC$ is extremely small, the stop band will be displayed.
\begin{figure}[htbp]
\centering
\scalebox{0.52}{\includegraphics{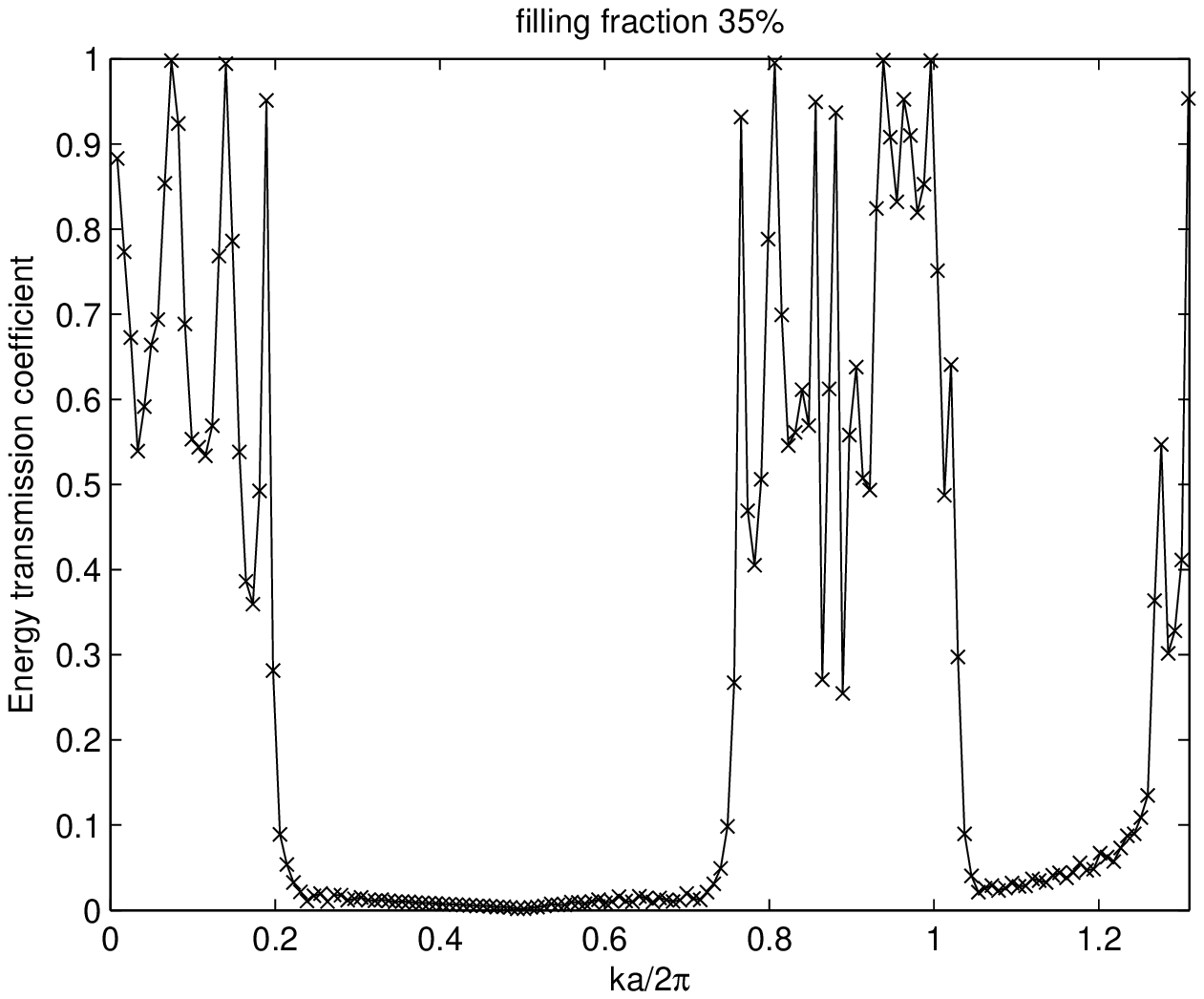}}\scalebox{0.52}{\includegraphics{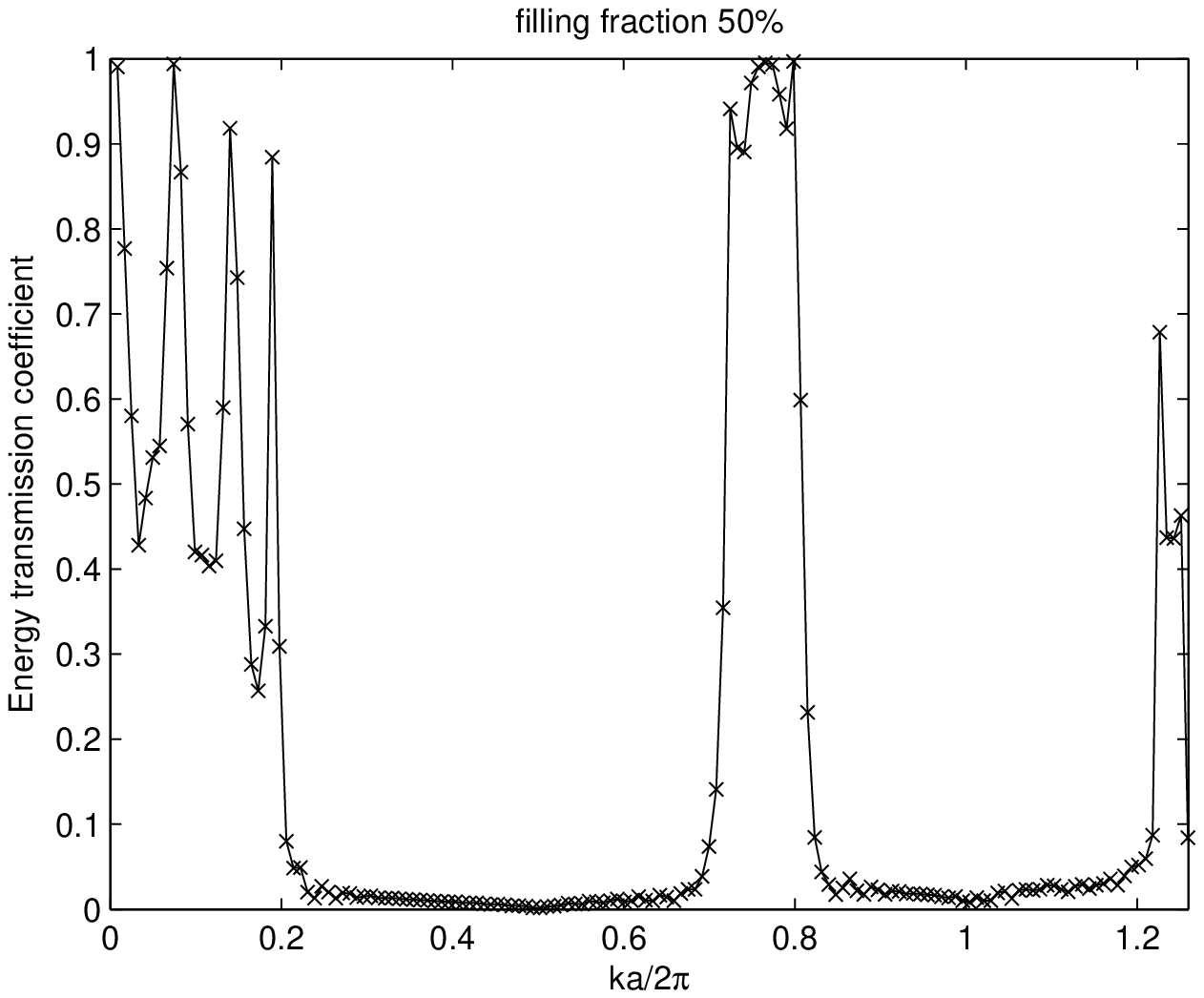}}\\
{\small Fig.6 $ETC$ for the water/mercury system with circular scatterer in square lattice.}
\end{figure}
\begin{figure}[htbp]
\centering
\scalebox{0.52}{\includegraphics{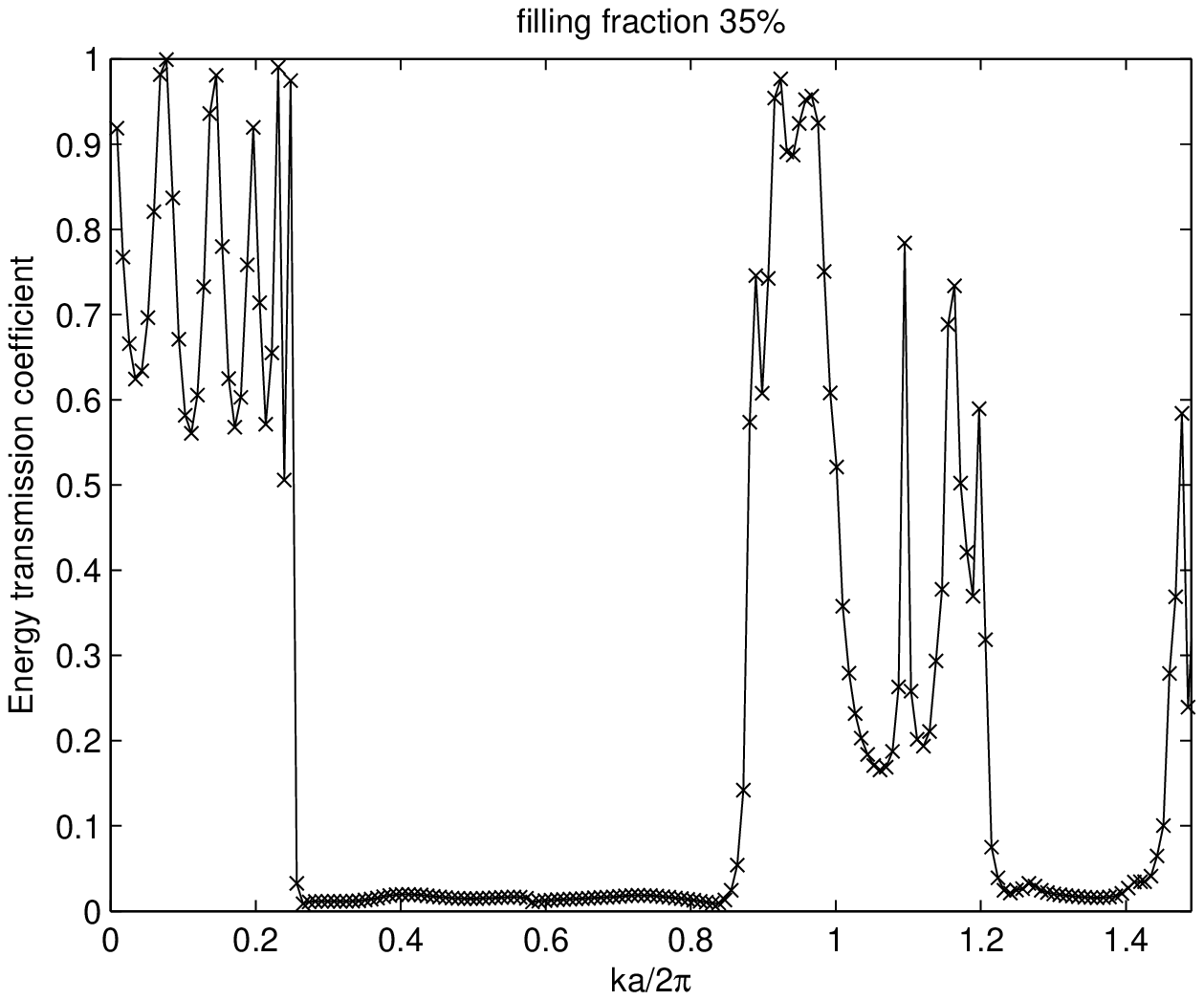}}\scalebox{0.52}{\includegraphics{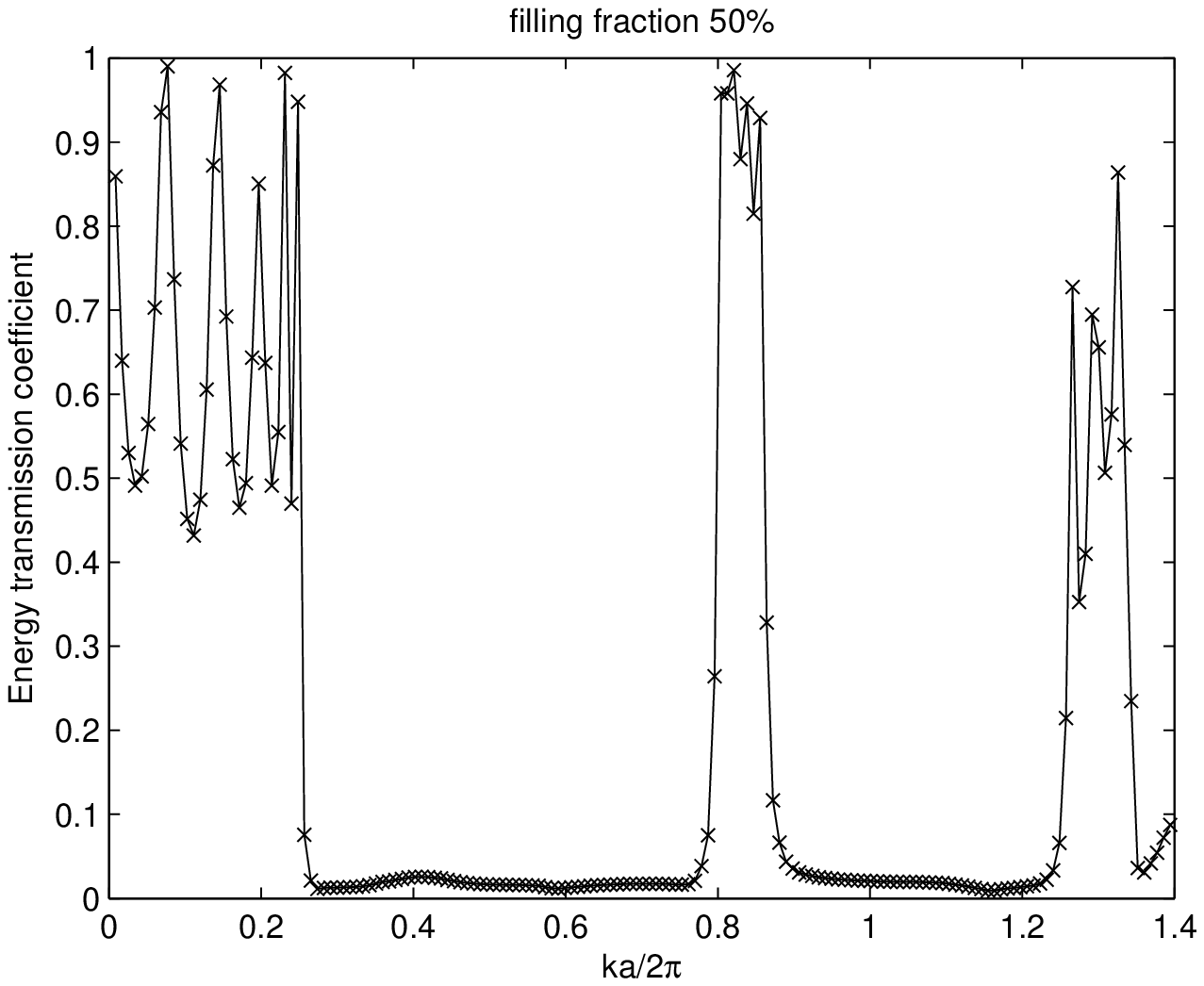}}\\
{\small Fig.7 $ETC$ for the water/mercury system with circular scatterer in hexagon lattice.}
\end{figure}
\begin{figure}[htbp]
\centering
\scalebox{0.52}{\includegraphics{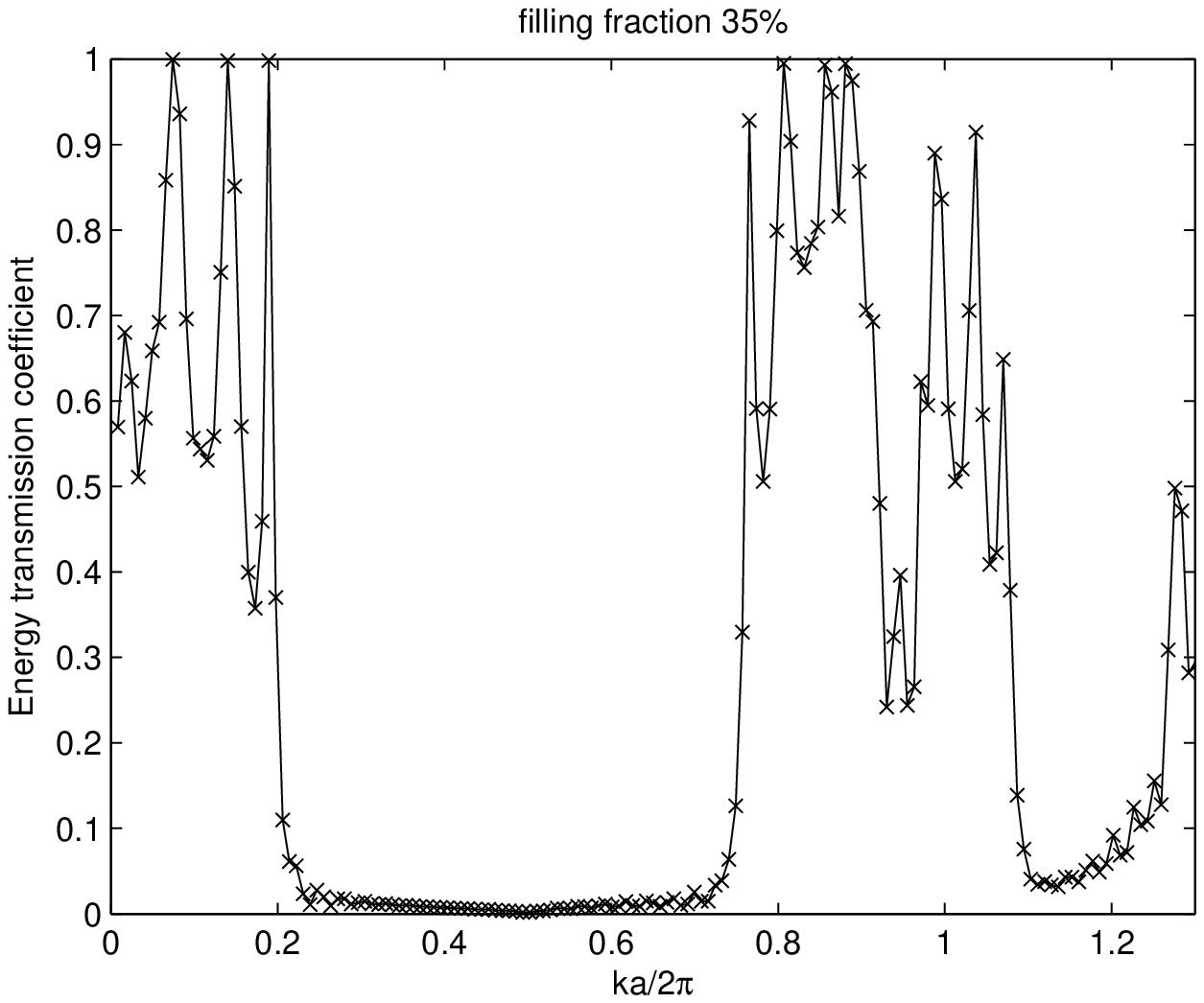}}\scalebox{0.52}{\includegraphics{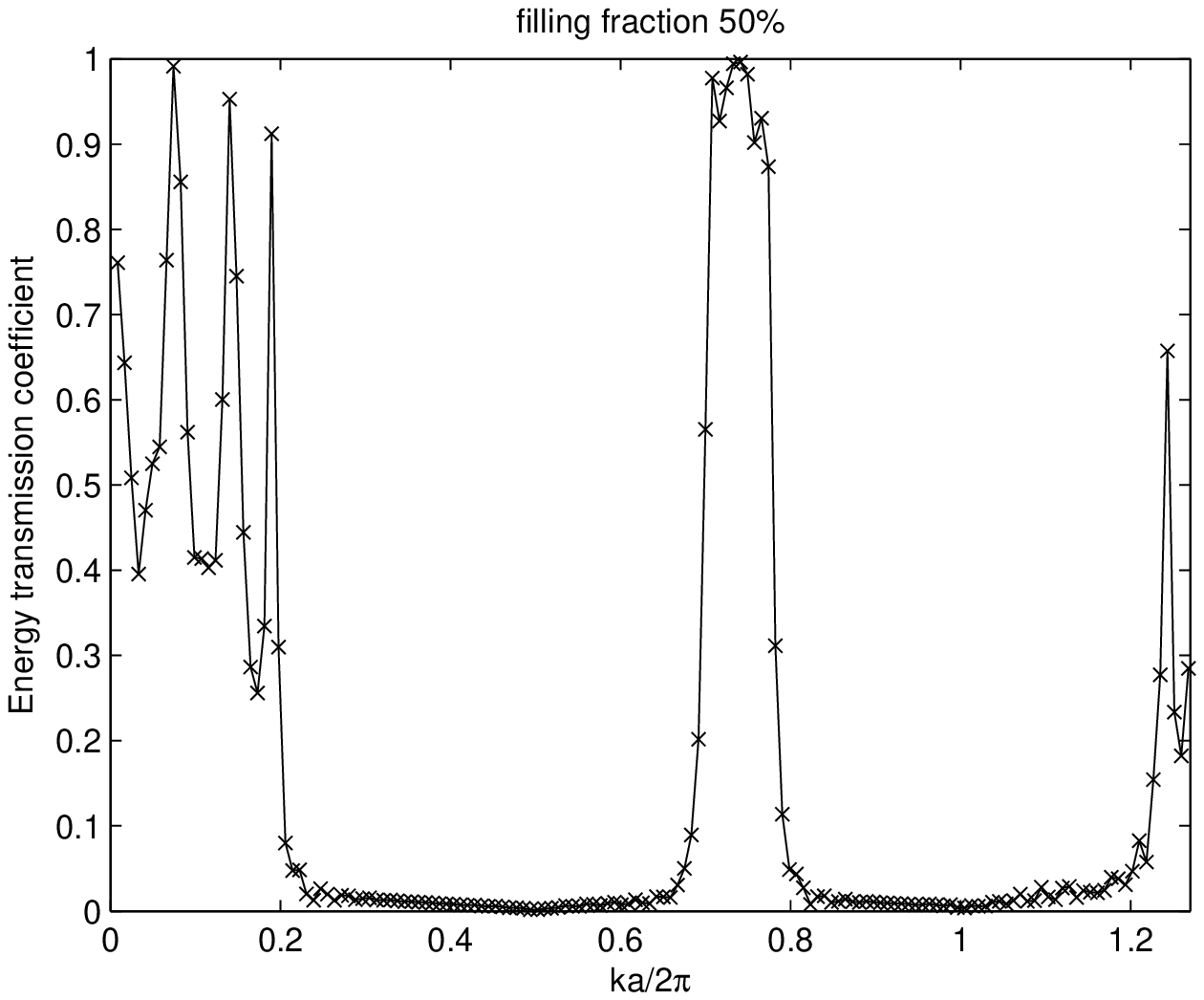}}\\
{\small Fig.8 $ETC$ for the water/mercury system with the round square scatterer in a square lattice.}
\end{figure}

The first numerical experiment is about the water/mercury system (the arrays of water scatterers are in the mercury matrix) with circular scatterer in square lattice, where the lattice constant $a=0.5$, $M_H=6, M_L=4$, the boundary of each scatterer is discretized to 256 points. Since the wave number $k=\omega/c$, it follows that $k=1.0345k_1$. In Fig.7, the energy transmission coefficients are plotted as the functions of $ka/2\pi$ for the filling fraction $35\%$ and $50\%$, the complete band gaps ($ETC<3\%$) are clearly shown in the figure.

The second numerical experiment is about the water/mercury system with circular scatterer in hexagon lattice, where the lattice constant is also $a=0.5$, and $M_H=6, M_L=4$. The results are plotted in Fig.8.

In the next numerical experiment, we select a round square scatterer with the boundary curve:
$$\Gamma: \big(r(3+2\sin^2\theta)\cos{\theta},r(3+2\cos^2\theta)\sin{\theta}\big),\qquad 0\leq\theta\leq2\pi,$$
the area of the scatterer is $23\pi{r}^2/2$. Fig.9 shows the $ETC$ for the water/mercury system with the round square scatterer in a square lattice.

The detailed band gaps for the above three numerical experiments are shown in Table 3. These results are very close to those obtained by PWE with 1681 plane waves \cite{17}.
\begin{center}
{Table 3. Band gaps for the water/mercury systems.}\\
{\begin{tabular}{l|c|c}\hline
 scatterer/lattice  & filling fraction $35\%$  & filling fraction $50\%$ \\\hline
 circular/square & ~$[0.2305,0.7244]$~ \vline\, ~$[1.0394,1.2178]$~ & ~$[0.2305,0.6915]$~ \vline\, ~$[0.8479,1.1772]$~ \\\hline
 circular/triangle & ~$[0.2652,0.8555]$~ \vline\, ~$[1.2319,1.4202]$~ & ~$[0.2567,0.7785]$~ \vline\, ~$[0.8983,1.2319]$~ \\\hline
 round square/square & ~$[0.2305,0.7244]$~ \vline\, ~$[1.0948,1.1607]$~ & ~$[0.2305,0.6668]$~ \vline\, ~$[0.8149,1.1689]$~ \\\hline
\end{tabular}}
\end{center}

\begin{figure}[htbp]
\centering
\scalebox{0.52}{\includegraphics{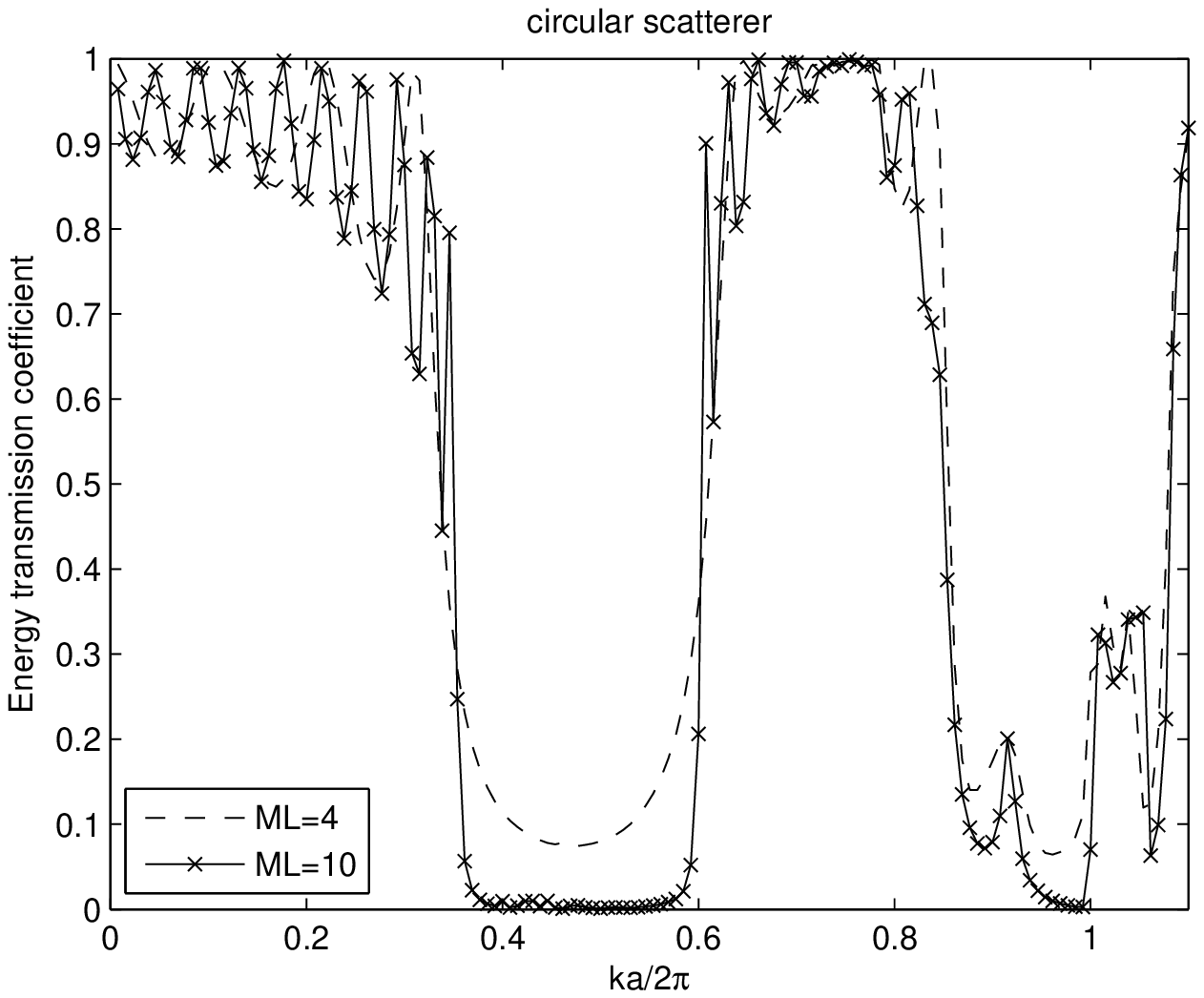}}\scalebox{0.52}{\includegraphics{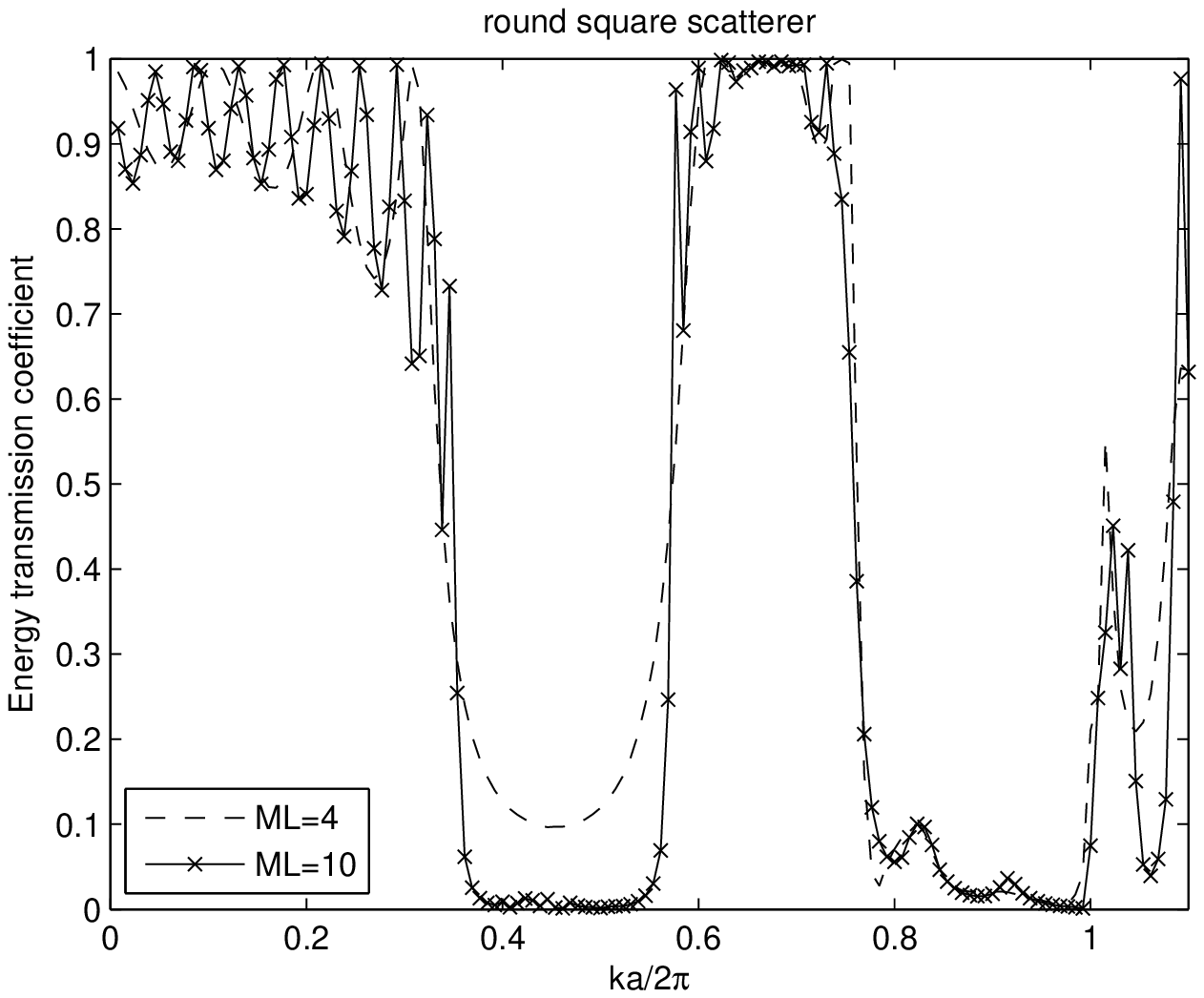}}\\
{\small Fig.9 $ETC$ for the mercury/water system with square lattice for filling fraction $35\%$.}
\end{figure}

At last, we consider the mercury/water system with circular scatterer in square lattice, the lattice constant is also $0.5$. The energy transmission coefficients are plotted in Fig.10. It is seen that the band gaps of mercury/water system are much narrower than that of water/mercury system. In addition, the stop band is not visible for $M_L=4$, with the increase of $M_L$ ($M_L=10$), not only the stop band is more visible, but also the edge of the stop band is sharper.

Numerical experiments show that the present FM-BEM can calculate the acoustic band gap of 2D liquid phononic crystal efficiently and accurately for arbitrary material combination and scatterers' shape. The idea of this paper can also be used to calculate the acoustic band gap of 2D solid phononic crystals and mixed solid-liquid phononic crystals.

It should be pointed out that, as can be seen from Fig.7-10, the results for high frequency seem to be unsatisfactory. This is because we choose the same truncation number $p$ for different frequencies. In fact, since $J_n(z)$ and $Y_n(z)$ are fluctuating functions when $n<z$, it follows that the truncation error of Graf's addition theorem is also fluctuating. With the increase of frequency, a larger truncation number $p$ should be taken to meet the accuracy requirement, this will also increase the computational complexity of the algorithm.

\section*{Acknowledgements}
This work is supported by the National Natural Science Foundation of China (11201373) and Natural Science Foundation of Shaanxi Provincial Department of Education (14JK1747).


\begin{thebibliography}{50}
\bibitem{1} D Colton,R Kress. Integral Equation Methods in Scattering Theory. John Wiley and sons, New York, 1983.
\bibitem{2} R Kress. Boundary integral equation in time-harmonic acoustic scattering. Math. Comput. Modelling, 15(3-5)(1991):229-243.
\bibitem{3} R Kress,G F Roach. Transmission Probelms for The Helmoltz Equation. J. Math. Phys., 19(1978):1433-1437.
\bibitem{4} R Kress. On the numerical solution of a hypersingular integral equation in scattering theory. J. Comput. Appl. Math., 61(1995):345-360.
\bibitem{5} V Rokhlin. Rapid solution of integral equations of classical potential theory. J. Comp. Phys., 60(1985):187-207.
\bibitem{6} V Rokhlin. Rapid solution of integral equations of scattering theory in two dimensions. J. Comp. Phys., 86(1990):414-439.
\bibitem{7} Y J Liu. Fast Multipole Boundary Element Method-Theory and Applications in Engineering. Cambridge University Press, Cambridge, 2009.
\bibitem{8} N Nishimura. Fast multipole accelerated boundary integral equation methods. Appl. Mech. Rev., 55(4)(2002):299-324.
\bibitem{9} S Amini,A T J Profit. Multi-level fast multipole solution of the scattering problem. Eng. Anal. Bound. Elem., 27(2003):547-564.
\bibitem{10} M Fischer,Ute Gauger,Lothar Gaul. A multipole Galerkin boundary element method for acoustics. Eng. Anal. Bound. Elem., 28(2004):155-162.
\bibitem{11} Kushwaha M S,Halevi P,G Mart¨ªnez,et al. Theory of acoustic band structure of periodic elastic composites. Phys. Rev. B, 1994,49(4):2313-2322.
\bibitem{12} F Kobayashi,S Biwa,N Ohno. Wave transmission characteristics in periodic media of finite length: multilayers and fiber arrays. Int. J. Solids Struct., 41(2004):7361-7375.
\bibitem{13} Tanaka Y,Tomoyasu Y,Tamura S I. Band structure of acoustic waves in phononic lattices: two-dimensional composites with large acoustic mismatch. Phys. Rev. B, 2000,62(11):7387-7392.
\bibitem{14} Kafesaki M, Economou E. Multiple-scattering theory for three-dimensional periodic acoustic composites. Phys. Rev. B, 17(1999):11993-12001.
\bibitem{15} P A Knipp, T L Reinecke. Boundary-element calculations of electromagnetic band-structure of photonic crystals. Phys. E, 2(1998):920-924.
\bibitem{16} A Barnett, L Greengard. A new integral representation for quasi-periodic fields and its application to two-dimensional band structure calculations. J. Comp. Phys., 229(2010):6898-6914.
\bibitem{17} A Barnett, L Greengard. A new integral representation for quasi-periodic scattering problems in two dimensions, BIT Numer. Math.,51(2011):67-90.
\bibitem{18} F L Li. Investigation on boundary integral equation method for calculation of band structures and transmission spectra of two-dimensional phononic crystals. Doctoral dissertation. Beijing Jiaotong University, Beijing, 2011.
\bibitem{19} Y Otani, N Nishimura. An FMM for periodic boundary value problems for cracks for Helmholtz' equation in 2D. Int. J. Numer. Meth. Engng, 73(2008):381-406.
\bibitem{20} Y Otani, N Nishimura. A periodic FMM for Maxwell's equations in 3D and its applications to problems related to photonic crystals. J. Comp. Phys., 227(2008):4630-4652.
\bibitem{21} Y Otani, N Nishimura. An FMM for orthotropic periodic boundary value problems for Maxwell¡¯s equations. Waves Random Complex, 19(2009):80¨C104.

\bibitem{22} W H Meng,L T Wang. Analysis of the convergence rates for the truncation errors of periodic Green's function of Helmholtz equations and its partial derivatives. Math. Numer. Sin., 37(2)(2015):123-136.
\bibitem{23} Frank W J Olver, NIST Handbook of Mathematical Functions, Cambridge University Press, 2010.
\bibitem{24} W H Meng. Bound on global error of the fast multipole method for Helmholtz equation in 2-D. arXiv:1806.08512 [math.NA].
\bibitem{25} N Yarvin and V Rokhlin. Generalized Gaussian quadratures and singular value decomposition of integral operators. SIAM J. Sci. Comput., 20(2)(1998):699-718.
\end{thebibliography}
\end{document}